\documentclass[11pt]{article}

\usepackage{natbib}
 \bibpunct[, ]{(}{)}{,}{a}{}{,}%

\usepackage{booktabs}
\usepackage{array}
\usepackage{fixltx2e}
\newcolumntype{P}[1]{>{\centering\arraybackslash}p{#1}}

\newcolumntype{C}[1]{>{\centering\let\newline\\\arraybackslash\hspace{0pt}}m{#1}}

\usepackage[a4paper,margin=0.6in]{geometry}
\usepackage{url,amsthm,amsmath,amsfonts,amssymb,amscd,yhmath,mathrsfs,algpseudocode,latexsym,enumerate,bm}
\usepackage{graphicx}
\usepackage{subfig}
\usepackage[font=small]{caption}

\usepackage{color,colortbl}
\usepackage{multirow}

\usepackage{symboldef}
\newtheorem{theorem}{Theorem}
\usepackage{todonotes}
\usepackage{rotating}
\date{}
\definecolor{Gray1}{gray}{0.8}
\definecolor{Gray1-5}{gray}{0.85}
\definecolor{Gray2}{gray}{0.9}

\begin{document}

\title{A computationally efficient Benders decomposition for energy systems planning problems with detailed operations and time-coupling constraints}

\author{Anna Jacobson\textsuperscript{\ref{LSI}} \\ 
annafj@princeton.edu \\ \\
Filippo Pecci \textsuperscript{\ref{Andlinger}} \\ \\
Nestor Sepulveda \textsuperscript{\ref{MIT},\ref{McKinsey}}\\ \\
Qingyu Xu \textsuperscript{\ref{EIRI}}\\ \\
Jesse Jenkins \textsuperscript{\ref{Andlinger},\ref{MAE}}}

\maketitle
\begin{enumerate}
    \item \label{LSI} Lewis-Sigler Institute of Integrative Genomics, Princeton University, Princeton, NJ, United States of America
    \item \label{Andlinger} Andlinger Center for Energy and the Environment, Princeton University, Princeton, NJ, United States of America
    \item \label{MIT} Research Affiliate, Massachusetts Institute of Technology, Cambridge, MA, Untied States of America
    \item \label{McKinsey} Manager, McKinsey and Company, United States of America
    \item \label{EIRI} Energy Internet Research Institute, Tsinghua University, Beijing, People’s Republic of China
    \item \label{MAE} Department of Mechanical and Aerospace Engineering, Princeton University, Princeton, NJ, United States of America
\end{enumerate}
\section*{Abstract}
Energy systems planning models identify least-cost strategies for expansion and operation of energy systems and provide decision support for investment, planning, regulation, and policy. Most are formulated as linear programming (LP) or mixed integer linear programming (MILP) problems. Despite the relative efficiency and maturity of LP and MILP solvers, large scale problems are often intractable without abstractions that impact quality of results and generalizability of findings. We consider a macro-energy systems planning problem with detailed operations and policy constraints and formulate a computationally efficient Benders decomposition separating investments from operations and decoupling operational timesteps using budgeting variables in the master model. This novel approach enables parallelization of operational subproblems and permits modeling of relevant constraints coupling decisions across time periods (e.g. policy constraints) within a decomposed framework. Runtime scales linearly with temporal resolution; tests demonstrate substantial runtime improvement for all MILP formulations and for some LP formulations depending on problem size relative to analagous monolithic models solved with state-of-the-art commercial solvers. Our algorithm is applicable to planning problems in other domains (e.g. water, transportation networks, production processes) and can solve large-scale problems otherwise intractable. We show that the increased resolution enabled by this algorithm mitigates structural uncertainty, improving recommendation accuracy.

\subsection*{Keywords}
macro-energy systems, capacity expansion planning, Benders decomposition, mixed integer linear programming, linear programming, decomposition methods

\section{Introduction}
Energy systems planning problems optimize resource investments, retirements, and operations to minimize total system cost (equivalently, maximize societal welfare) subject to technological, political, environmental, and economic constraints. These capacity expansion problems support decision-making in investment planning, regulation, and policy. Assuming perfect foresight and free entry, central-planning optimizations simulate a market under perfect competition; they are thus able to go beyond providing guidance on capacity deployment and retirement and play essential roles analyzing government policies~\citep{nzap,victoria2022,ricks2023} and advanced technologies’ role in decarbonized energy systems~\citep{victoria2020,mallapragada2020,ricks2022}. 

The vast majority of macro-energy systems planning problems (large-scale planning problems with regional or national scope) are implemented as linear programming (LP) or mixed integer linear programming (MILP) problems~\citep{ringkjob2018review,cho2022} due to the relative simplicity, computational efficiency, and maturity of LP and MILP solution methods and the fact that most salient system characteristics can be represented with reasonable accuracy using linear formulations.
Recently, the increasing penetration of variable renewable energy resources (VRE) has required much greater temporal, geospatial, and operational resolution to accurately capture key physical, economic, and engineering considerations that affect investment and retirement decisions~\citep{helisto2021impact}. Electrification of transportation, heating, and industrial processes and production of hydrogen from electrolysis is also increasing the relevance of electricity systems and more tightly-coupling electricity with other energy and industrial systems, requiring the formulation and solution of multi-sector energy systems models~\citep{brown2018,he2021}. As a result, full-resolution, full-scale models are composed of millions of variables and constraints, thereby risking intractability.

In order to deal with computational constraints, macro-energy systems planning problems are often heavily abstracted; most models either by downsample or subsample constituent time periods, aggregate regions into larger geographic zones, and/or relax operational constraints of physical systems. For a recent review, see \cite{cho2022}.

Simplified model structures can improve runtime but significantly impact investment and policy recommendations derived from energy planning models~\citep{bistline2022actions}. Modeling too few representative days may fail to capture weekly demand patterns and their influence on thermal plant unit commitment (UC) and storage dispatch decisions~\citep{mallapragada2020}. Employing non-sequential time slices may poorly incorporate weather patterns~\citep{poncelet2016impact} and prevent accurate modeling of flexibility requirements or energy storage. Temporal clustering is shown to significantly impact investment and operation of VRE~\citep{pfenninger2017dealing}. Simplification of system operation has been shown to affect investment recommendations in transmission planning~\citep{xu2019value,neumann2022} and dispatch decisions for systems with UC constraints~\citep{poncelet2020unit,palmintier2013heterogeneous}. Different means of geospatial aggregation techniques also impact model output~\citep{siala2019impact,frysztacki2022}. Because of these effects, abstractions must be carefully tailored to each study's focus; this inhibits the generalizability of any given model and the replicability of solutions when multiple models' solutions are compared.

Energy systems planning problems are often intractable even while deploying significant abstractions; these models require computationally efficient solution methods to terminate. \cite{lohmann2017tailored} developed tailored Benders decomposition algorithms for three simplified planning problems. Among the cases studied, two did not include time coupling constraints, as this omission allowed decomposed subproblems to be solved in parallel once investment decision variables had been fixed. The third case included just a single aggregated demand zone and omitted transmission operations and investments, storage resources, and policy constraints. Multi-period planning problems solved in \cite{lara2018deterministic} and \cite{li2022mixed} considered detailed operational and time-coupling constraints and decomposed the resulting MILP problems into a series of optimization problems per planning period. These studies did not investigate computationally efficient methods to solve the single-period subproblems; as a result, the largest case study in \cite{lara2018deterministic} and \cite{li2022mixed} had only $6$ zones and modeled each planning year using only $15$ representative days. The work by \cite{munoz2016new}, is one of the few to investigate computationally efficient methods for single-period planning problems. This study considered a full operational year with hourly resolution but ignored storage resources, ramping limits, and UC constraints, all of which couple operational decisions across time periods. \cite{sepulveda2020decarbonization} proposed a nested decomposition algorithm to solve a single-year planning problem with detailed operational and time coupling constraints. In the first decomposition stage, investment decisions are separated from operational decisions. Then, the operational subproblem is solved using a Dantzig-Wolfe decomposition. However, such a technique does not allow for cuts to be fully decoupled by timestep and requires each iteration of the outer decomposition algorithm to await convergence of the inner decomposition before continuing to the next iteration.

Our study goes further than state-of-the-art models in the literature by investigating decomposition methods for a single-period energy systems planning problem with hourly resolution and detailed operation and time-coupling constraints. We aim to mitigate the errors introduced by temporal clustering by solving the planning problem for a full operational year with minimal downscaling or subsampling of intra-annual timesteps and to minimize abstractions of operational constraints. The model described below is an electricity system planning problem with detailed operational decisions and constraints on ramping, storage operations, and start-up and shut-down~(UC) for thermal resources. Our formulation further includes policy constraints that couple time steps across the planning period, e.g. caps on annual CO$_2$ emissions or a renewable portfolio standard~(RPS). Decision variables consist of generation, energy storage and transmission investment and retirement decisions and operational decisions like resource dispatch. For MILP cases, we constrain all investment decisions to be integer. In comparison to previous literature~\citep{lara2018deterministic,li2022mixed}, we solve the resulting energy systems planning problem for a full operational year with hourly resolution, resulting in $8,736$ time steps (52 weeks). 

In this work, we propose a new decomposition scheme to separate investment decisions into a master model, along with a series of budgeting variables representing the time-coupling constraints (namely, CO$_2$ and RPS constraints.) This representation of policy constraints using budgeting variables is novel, to the best of the authors' knowledge. The reformulation described below enables the decomposition of a monolithic problem into a master model and several operational subproblems. Because each subproblem is fully decoupled, operational subproblems are solved in parallel with solutions incorporated as one Benders cut per subproblem per algorithmic iteration. As shown in Section~\ref{results}, inclusion of multiple cuts per iteration significantly improves the method's computational performance compared to both monolithic solution approaches and standard, single-cut, Benders decomposition implementations. The inclusion of the budgeting variables for time-coupling constraints is thus a key novel contribution of this study.

The capacity expansion problems considered here belong to the wider framework of integrated planning in infrastructure systems. These optimization problems are characterized by a complex structure wherein diagonal blocks of the matrix composing systems' variables and constraints are linked by both complicating variables (e.g., investments) and complicating constraints (e.g., policy constraints) see Figure~\ref{fig:probstruct}. Examples of other problems in this class include but are not limited to: optimal production planning in industrial processes, where inventory constraints link together different planning periods~\citep{shah2012}; stochastic planning of water and wastewater systems~\citep{naderi2017}; and optimal placement of electric bus charging stations, where both locations and scheduling are optimized~\citep{an2020}. In all of these cases, the resulting optimization problems are too large to be solved by monolithic approaches and tailored decomposition methods such as the one described here are required. Because our decomposition scheme is not restricted to the framework of energy systems, our computational experiments can inform researchers working on integrated planning problems in other application areas.

Section~\ref{sec:probform} describes the formulation of the energy systems planning problem with detailed operational and time-coupling constraints. In Section \ref{sec:solution_method}, we develop a novel Benders decomposition scheme which divides a full operational year into subperiods which can be processed in parallel. In Section \ref{sec:num_exp}, we evaluate the efficiency of the developed method using case studies derived from the Eastern United States with varying spatial extent ranging from $2$ - $19$ zones and with levels of temporal resolution ranging from $2$ to $52$ weeks. We conclude by examining the impacts of our enabled increased temporal resolution on resources' investment recommendations.

\section{Problem formulation}
\label{sec:probform}
We formulate an energy systems planning problem wherein electricity generation, energy storage, capacity expansion and retirement, and energy dispatch are jointly optimized over a single planning period. An in-depth description of the constraints comprising the problem are included in the appendix, sections \ref{prob_description_inv} - \ref{prob_description_policy}. The objective function being minimized represents investment and operational costs and includes penalties for violating policy constraints. In this work, we consider three different policy scenarios and their corresponding policy constraints, as described in Table~\ref{tab:policy_cases}.
 \begin{table}[h!]
 \centering
 \renewcommand{\arraystretch}{1.45}
 \footnotesize{
 \caption{Policy scenarios considered in this study.}
 \label{tab:policy_cases}
     \renewcommand{\arraystretch}{1.45}
    \begin{tabular}{cl}
    \toprule
    Scenario & Description \\
    \midrule
    REF & The reference case. Emissions and dispatch by resource type are unrestricted. \\
    RPS & Renewable portfolio standard. 70\% of generation must come from qualifying resources (e.g., VRE).\\
    CO$_2$ & CO$_2$ emissions cap. Emissions are constrained to 0.05 tons per MWh.\\
    \bottomrule
    \end{tabular}}
\end{table}

Optimization constraints (see section \ref{prob_description_op}) include maximum output limits for renewable resources, ramping limits, and start-up and shut-down restrictions for thermal generators (e.g. UC constraints.) As noted by~\cite{palmintier2015impact}, the ability to represent UC significantly impacts model results and improves accuracy. Here, we implement an aggregated UC model similar to \cite{palmintier2013heterogeneous,palmintier2015impact} and relax the integer constraints on UC variables to further improve scalability.

In most cases, planning problems for energy networks with several thousand nodes (i.e., buses) are intractable; we therefore divide the geographical area of interest into zones which incorporate real-world data (see section \ref{powergenome_description}) on demand profiles and climate conditions. For each zone, we create clusters for resources (e.g., generators or storage units) based on technology, cost of connection to transmission grids, and operating parameters. We assume that resources within a given cluster have the same capacity size per unit and use integer variables to model the decisions to install or retire a number of units within each cluster. We also include interregional power transmission and integer investment decision variables for capacity expansion of existing transmission paths between modeled zones. 
In this way, the energy system is modeled as a graph where each node represents a zone with demand and constituent resources (e.g., generators, storage) and edges represent interzonal transmission capacity. 

In this study, we implement a transport model and do not consider Kirchoff's Voltage Law (KVL) or power losses. As noted in \cite{neumann2022}, inclusion of these features would have greater impact on systems with lower levels of spatial aggregation (e.g., hundreds of nodes) than those described in Section~\ref{sec:num_exp}, which include up to $19$ aggregated demand zones. However, the novel Benders decomposition method developed in Section~\ref{sec:solution_method} can be applied to planning problems with KVL constraints due to the separation of transmission investments and operations without the need for nonlinear solution methods; this is the subject of future work.

Macro-energy systems planning problems are commonly formulated using a selection of operational subperiods due to models' size~\citep{frew2016temporal,mallapragada2018,lara2018deterministic,mallapragada2020,li2022mixed}. Storage resources' operation and thermal resource UC and ramping constraints are often modeled within each subperiod by linking first and last time steps in a method known as ``circular indexing." This approximation assumes that storage levels and UC decisions across two subperiods are decoupled. Errors may arise when subperiods are too short~\citep{mallapragada2020}, as it becomes impossible to fully capture weekly demand and weather patterns and their influence on UC and storage dispatch decisions. This is often the case in previous studies~\citep{mallapragada2018,lara2018deterministic,li2022mixed}, where only a few representative 24-hour periods (days) are used to model systems' yearly operation. For the work in this manuscript, subperiods are one week long; circular indexing occurs over a 168-hour timeperiod.

We develop a Benders decomposition scheme to solve an optimization problem for a planning period with up to $52$ consecutive weeks of operational decisions. The investigation of appropriate timeseries clustering methods to select representative subperiods is outside the scope of this manuscript. For problems with fewer than 52 weeks we assume that a clustering method has been implemented, resulting in sampled subperiods indexed by set $W$. For each subperiod $w \in W$, we define its set of hours as $H_w = \{(w-1)\hoursperperiod_w + 1,\ldots, w\hoursperperiod_w\}$, where $\hoursperperiod_w$ is the number of time steps within the subperiod (in our case, $\hoursperperiod_w=168$).

\subsection{Overall problem formulation}
Appendix section \ref{prob_description} contains a list of all constraints included in the model. To highlight problem structure, we introduce a compact formulation below. Assume that when $\bm{u}$ and $\bm{v}$ are vectors, the inequality $\bm{u} \leq \bm{v}$ is evaluated component-wise. Let $\bm{y} \in \mathbb{R}^{m}$ be a vector grouping all investment decision variables, and let $\bm{R}$ and $\bm{r}$ be such that constraints \eqref{eq:invest_cons_1}-\eqref{eq:invest_cons_3} correspond to $\bm{R}\bm{y} \leq \bm{r}$. In addition, let vector $\bm{c}_I$ be such that the fixed cost objective terms~\eqref{eq:fixed_cost} are denoted $\bm{c}_I^T\bm{y}$. For each subperiod $w \in W$, consider a vector $\bm{x}_w \in \mathbb{R}^{n}$ grouping all operational decision variables, and let matrices $\bm{A}_w$ and $\bm{B}_w$, and vector $\bm{b}_w$ be such that $\bm{A}_w\bm{x}_w + \bm{B}_w\bm{y} \leq \bm{b}_w$ corresponds to constraints \eqref{eq:dem_balance}-\eqref{eq:operation_cons_9}. Let vector $\bm{c}_w$ be such that the objective function terms \eqref{eq:var_cost} + \eqref{eq:startup_cost} + \eqref{eq:nse_cost} + \eqref{eq:policy_cost} are equal to $\sum_{w \in W}\bm{c}_w^T\bm{x}_w$. Finally, let matrix $\bm{Q}_w$ and vector $\bm{e}$ be such that $\sum_{w\in W}\bm{Q}_w\bm{x}_w \leq \bm{e}$ represents the policy constraints across the different scenarios, (corresponding to \eqref{eq:rps_cons} for case RPS, to \eqref{eq:co2_cons} for case CO$_2$, and remaining unenforced for case REF.) The resulting MILP problem is:
\begin{subequations}
\label{eq:prob}
\begin{alignat}{3}
    &\text{minimize}&\; \; &\bm{c}_I^T\bm{y} + \sum_{w \in W}\bm{c}_w^T\bm{x}_w \\
    &\text{subject to}& & \bm{A}_w\bm{x}_w + \bm{B}_w\bm{y}\leq \bm{b}_w, \quad \forall w \in W\\
    &&& \sum_{w \in W}\bm{Q}_w\bm{x}_w \leq \bm{e} \label{eq:prob_coupling} \\
    &&& \bm{R}\bm{y} \leq \bm{r} \label{eq:probinv}\\
    &&& \bm{x}_w \geq 0,  \quad \forall w \in W \\
    &&& \bm{y} \geq 0 \\
    &&& \bm{y}  \in \mathbb{Z}^{m}.
\end{alignat}
\end{subequations}

\section{Solution method}
\label{sec:solution_method}
\begin{figure}[h]
\centering
\includegraphics[width=0.5\textwidth]{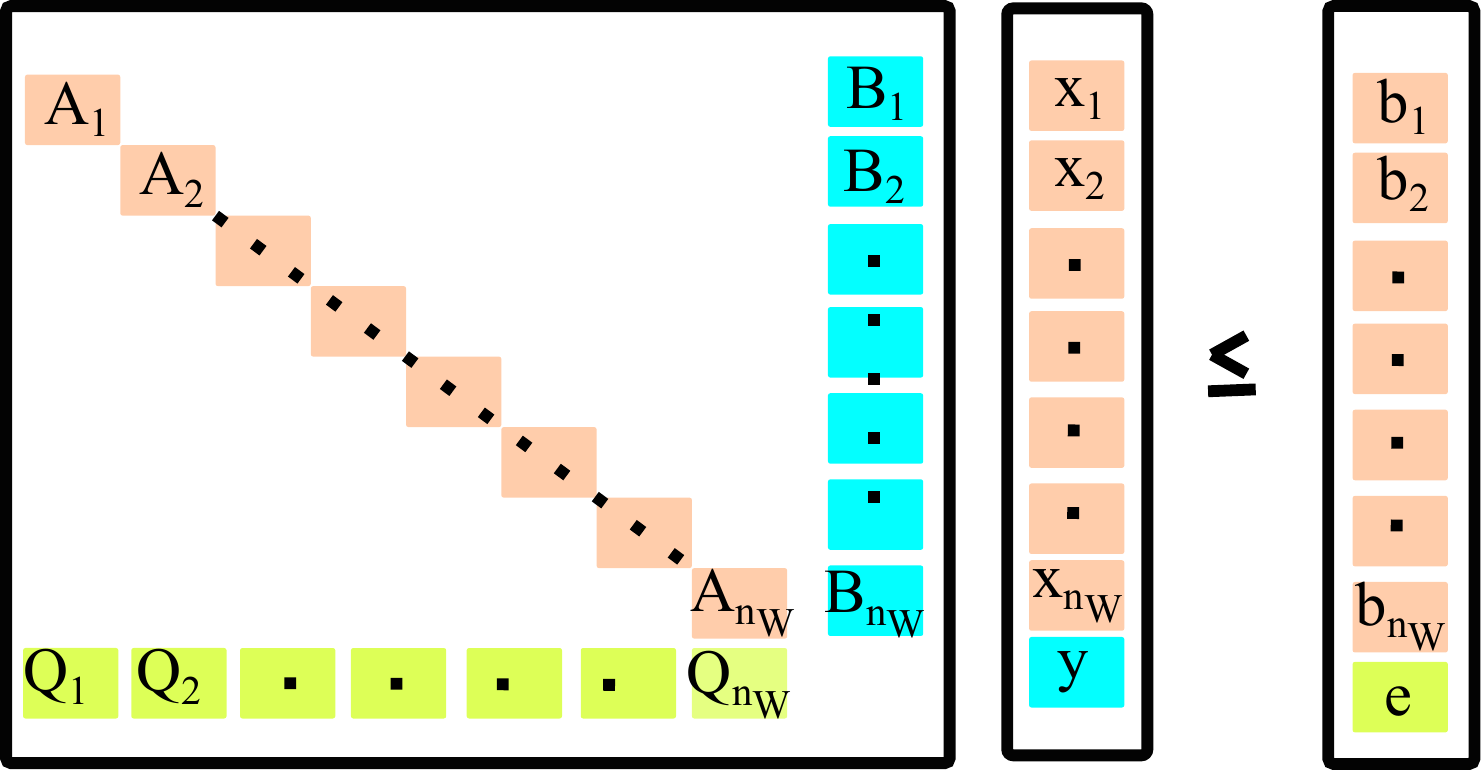}
\caption{Block structure of Problem~\eqref{eq:prob} with both complicating constraints and variables, where $n_W=|W|$. Investment-only constraints~\eqref{eq:probinv} are not pictured.}
\label{fig:probstruct}
\end{figure}

Problem~\eqref{eq:prob} is difficult to solve because it includes both complicating variables (e.g., investment decisions) and constraints (e.g., CO$_2$ limits) that link all operational subperiods - see Figure~\ref{fig:probstruct}. 

The complex structure of Problem~\eqref{eq:prob} is shared by several integrated planning problems in a variety of application areas~(e.g., \cite{shah2012,naderi2017,lara2018deterministic,an2020}). While ubiquitous, analogous optimization problems are particularly common in macro-energy systems modeling.
Multi-period planning problems in \cite{lara2018deterministic} and \cite{li2022mixed} are decomposed into a series of of single-period operational problems that are a special case of \eqref{eq:prob}, where investment decision variables are fixed. In these studies, the single-period operational problems were not further decomposed due to the presence of policy constraints like those seen in \eqref{eq:rps_cons} and \eqref{eq:co2_cons}. Few studies have developed efficient methods for single-period planning problems: \cite{munoz2016new} investigated the solution of a special case of Problem~\eqref{eq:prob} where storage, ramping limits, and thermal plant UC are not considered. However, the exclusion of UC as a system characteristic affects the accuracy of model results~\citep{palmintier2015impact}. Studies like the one performed by~\cite{lohmann2017tailored} included simplified cases of Problem~\eqref{eq:prob}, where either time-coupling or inter-zonal transmission constraints are ignored.

In the following, we develop a Benders decomposition scheme for Problem~\eqref{eq:prob} which enables parallel computation of subperiods. A standard implementation of Benders decomposition to Problem~\eqref{eq:prob} would consider only investment decisions $\bm{y}$ as complicating variables, and rely on the solution of a large, monolithic operational subproblem, spanning the whole year, at each iteration. This is the same approach used in \cite{lara2018deterministic} and \cite{li2022mixed} to solve analogous energy systems planning problems. Nested decompositions have also been proposed~\citep{sepulveda2020decarbonization} as means of further decomposing operational subproblems, but this technique still takes the operational problem as a single entity in the scope of the master model. In both cases, the master problem is given by:
\begin{subequations}
\label{eq:master_classic}
\begin{alignat}{3}
    &\text{minimize}&\; \; &\bm{c}_I^T\bm{y} + \theta\\
    &\text{subject to}& & \theta \geq \sum_{w \in W}f^{j}_w + (\bm{\pi}^j)^T(\bm{y}-\bm{y}^j), \quad \forall j=0,\ldots,k-1,  \label{eq:master_classic_cuts}\\
    &&& \bm{R}\bm{y} \leq \bm{r} \\
    &&& \bm{y} \geq 0 \\
    &&& \bm{y}  \in \mathbb{Z}^{m}.
\end{alignat}
\end{subequations}
Where $\theta$ represents an approximation of the operational cost within the master model. Benders cuts~\eqref{eq:master_classic_cuts} approximate operational cost based on a single operational subproblem, which models dispatch across the entire timeseries:
\begin{subequations}
\label{eq:subprob_classic}
\begin{alignat}{3}
    &\text{minimize}&\; \; &\sum_{w \in W}\bm{c}_w^T\bm{x}_w \\
    &\text{subject to}& & \bm{A}_w\bm{x}_w + \bm{B}_w\bm{y}\leq \bm{b}_w, \quad \forall w \in W \\
    &&& \sum_{w \in W}\bm{Q}_w\bm{x}_w \leq \mathbf{e} \label{eq:subprob_classic_coupling}\\
    &&& \bm{x}_w \geq 0,  \quad \forall w \in W \\
    &&& \bm{y} = \bm{y}^{k}  \quad \quad  :\bm{\pi}
\end{alignat}
\end{subequations}
Problem~\eqref{eq:subprob_classic} is not separable with respect to the subperiod index $w \in W$ due to constraints in~\eqref{eq:subprob_classic_coupling} that tie together all subperiods. As shown in Section~\ref{sec:num_exp}, modeling operations as a monolithic problem is impractical as Problem~\eqref{eq:subprob_classic} quickly becomes intractable as number of zones and subperiods increases. \cite{sepulveda2020decarbonization} suggested solving Problem~\eqref{eq:subprob_classic} using a nested Dantzig-Wolfe decomposition algorithm. However, such a strategy requires each iteration of the Benders decomposition algorithm to await convergence of the inner Dantzig-Wolfe decomposition before continuing. It also does not allow the inclusion of multiple decoupled cuts~\eqref{eq:master_classic_cuts} per iteration, resulting in increased number of iterations as seen in Table~\ref{fullop_table}.

In the following we take a different approach, enabling the separation of Problem~\eqref{eq:prob} into subproblems that can be solved in parallel without requiring a nested decomposition.  First, we prove the following result.
\begin{theorem}
Problem~\eqref{eq:prob} is equivalent to:
\begin{equation}
\label{eq:prob_decomp}
\begin{alignedat}{3}
    &\text{minimize}&\; \; &\bm{c}_I^T\bm{y} + \sum_{w \in W}\bm{c}_w^T\bm{x}_w \\
    &\text{subject to}& & \bm{A}_w\bm{x}_w + \bm{B}_w\bm{y}\leq \bm{b}_w, \quad \forall w \in W\\
    &&& \bm{Q}_w\bm{x}_w \leq \bm{q}_w, \quad w \in W \\
    &&& \sum_{w \in W}\bm{q}_w = \bm{e}\\
    &&& \bm{R}\bm{y} \leq \bm{r} \\
    &&& \bm{x}_w \geq 0,  \quad \forall w \in W \\
    &&& \bm{y} \geq 0 \\
    &&& \bm{y}  \in \mathbb{Z}^{m}.
\end{alignedat}
\end{equation}
\end{theorem}
\proof{Proof.} 
Without loss of generality, we assume that $W = \{1,\ldots,n_W\}$. We will show that:
\begin{equation}
\label{eq:main_result}
\sum_{i=1}^{n_W}\bm{Q}_i\bm{x}_i \leq \bm{e} \Leftrightarrow  \exists (\bm{q}_i)_{i=1}^{n_W} \text{ such that } \sum_{i=1}^{n_W}\bm{q}_i=\bm{e} \text{ and } \bm{Q}_i \bm{x}_i \leq \bm{q}_i,\quad \forall i =1,\ldots,n_W.
\end{equation}
The implication $``\Leftarrow"$ follows by definition of vectors $\bm{q}_1,\ldots,\bm{q}_{n_W}$. By induction, we show that also the $``\Rightarrow"$ implication holds. \\
\noindent Assume $n_W=2$ and $\bm{Q}_1\bm{x}_1 + \bm{Q}_2\bm{x}_2 \leq \bm{e}$. Then, if we define $\bm{q}_1=\bm{e}-\bm{Q}_2\bm{x}_2$ and $\bm{q}_2=\bm{Q}_2\bm{x}_2$ we have
\begin{equation}
\begin{split}
&\bm{q}_1+\bm{q}_2 = \bm{e} \\
&\bm{Q}_1\bm{x}_1 \leq \bm{e}-\bm{Q}_2\bm{x}_2 = \bm{q}_1\\
&\bm{Q}_2\bm{x}_2 =\bm{q}_2
\end{split}
\end{equation}
Next, we assume that it holds for $n_W=l$ and prove that this implies that it holds for $n_W=l+1$. We want to show that:
\begin{equation}
\sum_{i=1}^{l+1}\bm{Q}_i\bm{x}_i \leq \bm{e} \Rightarrow \exists \bm{q}_1,\ldots,\bm{q}_{l+1}\text{ such that } \sum_{i=1}^{l+1}\bm{q}_i=\bm{e}\text{ and } \bm{Q}_i\bm{x}_i\leq \bm{q}_i, \; \forall i=1\ldots,l+1
\end{equation}
Define $\bm{q}_{l+1}=\bm{Q}_{l+1}\bm{x}_{l+1}$ and observe: 
\begin{equation}
\sum_{i=1}^{l}\bm{Q}_i\bm{x}_i \leq \bm{e}-\bm{Q}_{l+1}\bm{x}_{l+1}=\bm{e}-\bm{q}_{l+1}.
\end{equation}
Since we are assuming that the claim holds for $n_W=l$, we have that there exist $\bm{q}_1,\ldots,\bm{q}_l$ such that $\sum_{i=1}^l\bm{q}_i = \bm{e}-\bm{q}_{l+1}$ and $\bm{Q}_i\bm{x}_i\leq \bm{q}_i$ for all $i=1,\ldots,l$, which completes the induction step.
\endproof
\subsection{Benders decomposition algorithm}
The structure of Problem~\eqref{eq:prob_decomp} suggests the implementation of a Benders decomposition algorithm wherein budgeting variables ($\bm{q}_w$) are used to implement~\eqref{eq:prob_coupling} and both $\bm{y}$ and $\bm{q}_w$ are treated as complicating variables. We note that this budgeting approach should be extensible to other constraints coupling subperiods, such as long-duration storage state-of-charge or hydro reservoir levels or similar time-coupling inventory constraints, although these constraints are not implemented in the present numerical tests. Initialize $\text{UB}=\infty, f^0_w=0$, $\bm{\pi}^0=\bm{0}$,  $\bm{y}^0=\bm{0}$, $\bm{\lambda}_w^0=\bm{0}$, and $\bm{q}_w^0=\bm{0}$, for all $w \in W$. For all $k=1,\ldots K_{\max}$ proceed as follows:

\noindent \textbf{Step 1.} Obtain estimated solutions $\bm{y}^k$ and $\bm{q}_w^k$, for all $w \in W$ by solving the master problem:
\begin{subequations}
\label{eq:master}
\begin{alignat}{3}
    &\text{minimize}&\; \; &\bm{c}_I^T\bm{y} + \sum_{w \in W} \theta_w\\
    &\text{subject to}& & \theta_w \geq f^{j}_w + (\bm{\pi}^j)^T(\bm{y}-\bm{y}^j) +  (\bm{\lambda}_w^j)^T(\bm{q}_w - \bm{q}_w^j), \quad \forall j=0,\ldots,k-1, \; w \in W \ \label{eq:master_cut} \\
    &&& \sum_{w \in W}\bm{q}_w = \bm{e}\\
    &&& \bm{R}\bm{y} \leq \bm{r} \\
    &&& \bm{y} \geq 0 \\
    &&& \bm{y}  \in \mathbb{Z}^{m},
\end{alignat}
\end{subequations}
where $f^{j}_w$ represents the current operational cost given a fixed set of investments and $\bm{\pi}^j$ and $\bm{\lambda}_w^j$ represent the Lagrangian multipliers associated with investment decisions and policy constraints (e.g., RPS) respectively. Set $\text{LB}$ to be the optimal value of Problem~\eqref{eq:master}.

\noindent \textbf{Step 2.} For every $w \in W$, solve the following linear subproblem:
\begin{equation}
\label{eq:subprob}
\begin{alignedat}{3}
    &\text{minimize}&\; \; &\bm{c}_w^T\bm{x}_w \\
    &\text{subject to}& & \bm{A}_w\bm{x}_w + \bm{B}_w\bm{y}\leq \bm{b}_w \\
    &&& \bm{Q}_w\bm{x}_w \leq \bm{q}_w \\
    &&& \bm{x}_w \geq 0 \\
    &&& \bm{y} = \bm{y}^{k}  &: \bm{\pi}\\
    &&& \bm{q}_w = \bm{q}_w^{k} &: \bm{\lambda},
\end{alignedat}
\end{equation}
and compute the optimal value $f^k_w$ and Lagrangian multipliers $\bm{\pi}^k$ and $\bm{\lambda}_w^k$ given the fixed set of investment decisions ($\bm{y}^{k}$) and coupling constraint budgets ($\bm{q}^{k}_w$) which are constants in the scope of the subproblem. Set $\text{UB} = \min(\text{UB}_{\text{prev}},\bm{c}_I^T\bm{y}^{k} + \sum_{w \in W}f^k_w)$, where $\text{UB}_{\text{prev}}$ is the upper bound from the previous iteration of the algorithm, ($\infty$ if $k=1$.) If $\frac{\text{UB}-\text{LB}}{\text{LB}} \leq \texttt{Rel}_{\texttt{tol}}$, then stop. Else, set $k=k+1$ and go back to Step 1.

Observe that subproblems in Step 2 are separable and can be solved in parallel, solving $|W|$ smaller operational subproblems, each formulated over $\hoursperperiod_w$ timesteps. This comes at the cost of adding $|W|$ constraints (Benders cuts) to the master problem per iteration, thereby increasing the size of the master problem compared to a standard Benders implementation. Numerical experiments reported in the next section suggest that this trade-off is worthwhile, as the total computational effort is dominated by the solution of the operational subproblems rather than the master problem. Furthermore, incorporation of multiple cuts allows more information to be communicated to the master problem at each iteration of the algorithm, decreasing the number of iterations needed for convergence.
\section{Numerical experiments}
\label{sec:num_exp}
\begin{figure}[h!]
    \caption{IPM regions used in our numerical experiments.}
    \label{geography}
    \subfloat[Zones, 2-zone simulation]{\includegraphics[width=0.48\textwidth]{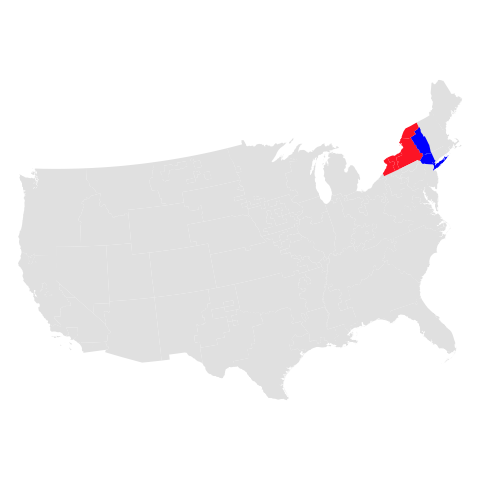}}
    \subfloat[Zones, 6-zone simulation]{\includegraphics[width=0.48\textwidth]{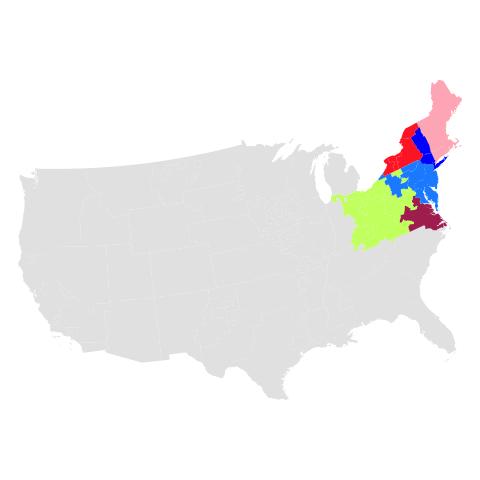}} \\
    \subfloat[Zones, 12-zone simulation]{\includegraphics[width=0.48\textwidth]{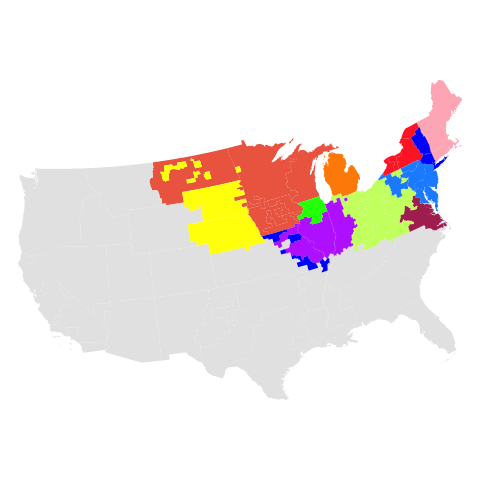}}
    \subfloat[Zones, 19-zone simulation]{\includegraphics[width=0.48\textwidth]{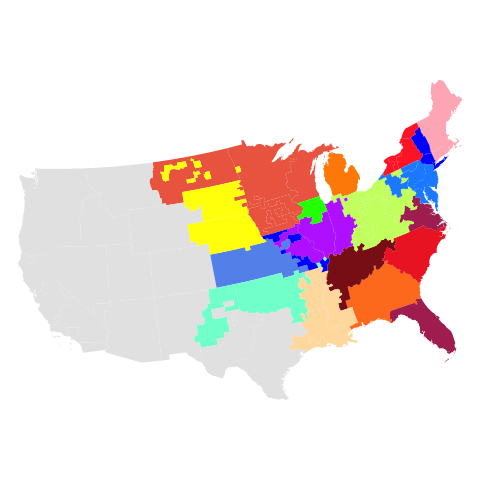}}
    \end{figure}

We consider 2-, 6-, 12-, and 19-zone cases of the Eastern United States, see Figure \ref{geography}. Initial capacity estimates are given for the year 2050 as output by the case generation software PowerGenome~\citep{schivley2021powergenome}, see section~\ref{powergenome_description}.
For each set of Eastern Interconnection regions, we represent the operational year using 2, 12, 22, 32, 42, and 52 representative weeks. Each planning problem is considered under the three policy scenarios described in Table \ref{tab:policy_cases}. The size of the resulting optimization problems is summarized in Table \ref{tab:problem_size}. 
\begin{table}[h!]
\renewcommand{\arraystretch}{1.45}
\caption{Model size by zone. Shows number of generator clusters, total number of variables and total number of constraints. Numbers of variables and constraints are those associated with the 52-week (full year) monolithic problem.}
\label{tab:problem_size}
\renewcommand{\arraystretch}{1.45}
\centering
\begin{tabular}{ccccc}
    \toprule
    Zones & $\vert G \vert$ & $\vert \UC \vert$ & Variables ($\cdot 10^9$) & Constraints ($\cdot 10^9$)\\
    \midrule
    2 & 62 & 16 & 1.1 & 3.4 \\
    6 & 175 & 54 & 3.4 & 10.5 \\
    12 & 285 & 106 & 6.2 & 19.3 \\
    19 & 437 & 167 & 9.7 & 30.4 \\
    \bottomrule
    \end{tabular}
\end{table}

In order to evaluate the impact of integer investment decision variables on the computational effort needed to solve energy systems planning problems, we consider two forms of Problem~\eqref{eq:prob}: in the first (hereafter ``MILP") investment and retirement variables are integer for both generation and transmission resources - this is the form presented in Problem~\eqref{eq:prob}. In the second (hereafter ``LP") we relax this constraint and allow all variables to be continuous. We compare the computational performance of the decomposition algorithm with a monolithic approach, where Problem~\eqref{eq:prob} is solved directly by state-of-the-art optimization solvers, and with a standard Benders decomposition algorithm where a single full operational subproblem is solved at each iteration.
\subsection{Computational setup}
All optimization problems are implemented in Julia 1.6.1~\citep{julia_2017}, where optimization solvers are called through JuMP 0.21.9~\citep{jump_2017}. We solve all MILP and LP problems using Gurobi (v9.0.1).
The simulations are run on Princeton University's Della computer cluster with 2.8 GHz Cascade Lake processors, Intel Xeon Platinum 8260 at 2.40 GHz. Problems are considered intractable if they require more than $48$ hours of computations or $200$ GB of memory to terminate.

We set a tolerance $\texttt{Rel}_{\texttt{tol}} = 10^{-3}$. Accordingly, optimality tolerance, MIP gaps, and barrier convergence tolerance for the monolithic models are set to $10^{-3}$. For these models, Gurobi is run using the barrier method with crossover turned off. This allows a fair comparison between the runtimes of monolithic and decomposed approaches, as Benders decomposition is not guaranteed to return basic solutions, and the crossover stage of the interior point solution method implemented by Gurobi requires considerable computational time to convert an optimal feasible solution to a an optimal basic solution. When solving the subproblem~\eqref{eq:subprob}, however, we enable crossover, as cuts~\eqref{eq:master_cut} require basic solutions.

We solve the operational subproblems~\eqref{eq:subprob} in parallel. For each iteration of the algorithm, $\vert W \vert$ CPUs are assigned one operational subproblem each. Data pertaining to cuts in~\eqref{eq:master_cut} is returned to the master model from these auxiliary CPUs. Returning only information necessary to formulate cuts (as opposed to entire solutions) allows minimal interprocess communication, thus speeding up time between iterations. Problem~\eqref{eq:master} waits for all subproblems to terminate before incorporating each of their cuts into a new iteration of the master problem and sending solutions $\bm{y}^{k}$ and $\bm{q}^{k}_w$ back to the subproblems for the next algorithmic iteration. For these cases, all subproblems were solved on a single compute node. In cases where the number of subproblems exceeds available CPUs on a single compute node, multiple subproblems may be run in series on each CPU, or multiple nodes could be coordinated via distributed parallelization (which incurs greater computational overhead for inter-process coordination across nodes).

\subsection{Results} \label{results}
We define runtime as the time spent initializing the problem, performing all computations within solvers, and writing solutions to file at the end of simulations. Table~\ref{runtime_table_co2} compares the computational performance of the Benders decomposition algorithm and the monolithic solution approach.

We find that using Benders decomposition brings little benefit for cases with small numbers of zones, especially for LP problems without integer variables. We observe a substantial reduction in computational time when considering MILP problems, particularly large scale problems. Additionally, all instances of the Bender's algorithm solved within the time limit of 48 hours whereas the monolithic algorithm saw 29 intractable MILP cases out of the 70 that were run. 

We observe that most intractable problem instances for the monolithic model occur in scenarios RPS and CO$_2$. Including these time-coupling constraints reflects common clean energy and emissions reductions targets; in our trials, the reference, RPS-, and CO$_2$-constrained cases incurred 8.6e8, 6.0e8, and 1.9e8 tons of CO$_2$ respectively across the year for the 19-zone case when simulated on 12 weeks of data. An inability to include these time-coupling constraints thus threatens models' relevance to systems undergoing decarbonization.
\begin{table}[h!]
\centering
 \renewcommand{\arraystretch}{1.45}

\caption{Runtime for Benders and monolithic models (100s) followed by ratio between monolithic and decomposition solution approaches. $\infty$ denotes a case with an intractable monolithic model. Cases that are intractable due to memory are noted with the superscript $M$. Cases that are intractable due to insufficient time are noted with a superscript $T$. Cases where the model outperformed its analogous model formulation are shown for the runtime rows, cases where the decomposed model outperformed monolithic are bolded for the ratio cases. Results shown for the CO$_2$-constrained case. Additional cases are included in the appendix.\label{runtime_table_co2}}
\footnotesize{
\begin{tabular}{c|r|r|P{0.045\textwidth}P{0.045\textwidth}P{0.045\textwidth}P{0.045\textwidth}P{0.045\textwidth}P{0.045\textwidth}|P{0.045\textwidth}P{0.045\textwidth}P{0.045\textwidth}P{0.045\textwidth}P{0.045\textwidth}P{0.045\textwidth}}
\toprule
\multicolumn{1}{c}{} & \multicolumn{1}{c}{$\vert Z \vert$} & \multicolumn{1}{c}{$\vert G \vert$} & \multicolumn{6}{c}{LP} & \multicolumn{6}{c}{MILP} \\
\cmidrule(lr){1-3} \cmidrule(lr){4-9} \cmidrule(lr){10-15}
\multicolumn{3}{r|}{Weeks $\rightarrow$} & 2 & 12 & 22 & 32 & 42 & 52 & 2 & 12 & 22 & 32 & 42 & 52 \\
\midrule
& 2 & 62 & 1.1 & 1.2 & \textbf{2.1} & \textbf{3.4} & \textbf{4.1} & \textbf{6.4} & \textbf{1.1} & \textbf{1.2} & \textbf{1.4} & \textbf{1.6} & \textbf{1.9} & \textbf{2.0} \\
& 6 & 175 & 5.9 & 10.1 & 16.6 & 19.7 & \textbf{26.0} & \textbf{28.3} & 6.1 & \textbf{10.6} & \textbf{16.4} & \textbf{21.6} & \textbf{25.9} & \textbf{36.8} \\
&12 & 285 & 66.9 & 79.3 & 112.8 & 135.1 & 151.6 & 173.0 & 89.5 & \textbf{86.2} & \textbf{128.0} & \textbf{165.2} & \textbf{191.1} & \textbf{200.4} \\
\parbox[t]{2mm}{\multirow{-4}{*}{\rotatebox[origin=c]{90}{Benders (100s)}}}$\; \;$&19 & 437 & 474.0 & 465.0 & 558.3 & 702.0 & 776.5 & 720.0 & 407.5 & \textbf{652.4} & \textbf{718.9} & \textbf{953.3} & \textbf{994.7} & \textbf{1123.4} \\
\midrule
\cellcolor{white}& 2 & 62 & \textbf{0.4} & \textbf{1.0} & 2.1 & 3.5 & 4.1 & 6.0 & 2.9 & 58.6 & 112.5 & 702.8 & 1344.3 & 651.9 \\
& 6 & 175 & \textbf{0.6} & \textbf{4.3} & \textbf{10.8} & \textbf{18.3} & 27.5 & 36.7 & 7.5 & 129.4 & $\infty^T$ & $\infty^M$ & $\infty^T$ & $\infty^T$ \\
& 12 & 285 & \textbf{1.2} & \textbf{10.7} & \textbf{27.9} & \textbf{50.4} & \textbf{79.3} & \textbf{106.8} & \textbf{24.5} & $\infty^T$ & $\infty^T$ & $\infty^T$ & $\infty^T$ & $\infty^T$ \\
\parbox[t]{2mm}{\multirow{-4}{*}{\rotatebox[origin=c]{90}{Mono. (100s)}}}$\; \;$&19 & 437 & \textbf{1.7} & \textbf{21.7} & \textbf{61.9} & \textbf{112.0} & \textbf{167.9} & \textbf{227.4} & \textbf{252.5} & $\infty^T$ & $\infty^T$ & $\infty^T$ & $\infty^T$ & $\infty^T$ \\
\midrule
& 2 & 62 & 0.3 & 0.8 & \textbf{1.5} & \textbf{2.2} & \textbf{2.3} & \textbf{3.2} & \textbf{2.5} & \textbf{48.0} & \textbf{78.2} & \textbf{444.8} & \textbf{722.7} & \textbf{329.2} \\
&6 & 175 & 0.1 & 0.4 & 0.7 & 0.9 & \textbf{1.1} & \textbf{1.3} & \textbf{1.2} & \textbf{12.3} & $\bm{\infty^T}$ & $\bm{\infty^M}$ & $\bm{\infty^T}$ & $\bm{\infty^T}$ \\
&12 & 285 & $<$0.1 & 0.1 & 0.3 & 0.4 & 0.5 & 0.6 & 0.3 & $\bm{\infty^T}$ & $\bm{\infty^T}$ & $\bm{\infty^T}$ & $\bm{\infty^T}$ & $\bm{\infty^T}$ \\
\parbox[t]{2mm}{\multirow{-4}{*}{\rotatebox[origin=c]{90}{Ratio}}}$\; \;$&19 & 437 & $<$0.1 & 0.1 & 0.1 & 0.2 & 0.2 & 0.3 & 0.6 & $\bm{\infty^T}$ & $\bm{\infty^T}$ & $\bm{\infty^T}$ & $\bm{\infty^T}$ & $\bm{\infty^T}$ \\
\bottomrule
\end{tabular}
}
\end{table}

\begin{figure}[h!]
\centering
    \includegraphics[width=0.8\textwidth]{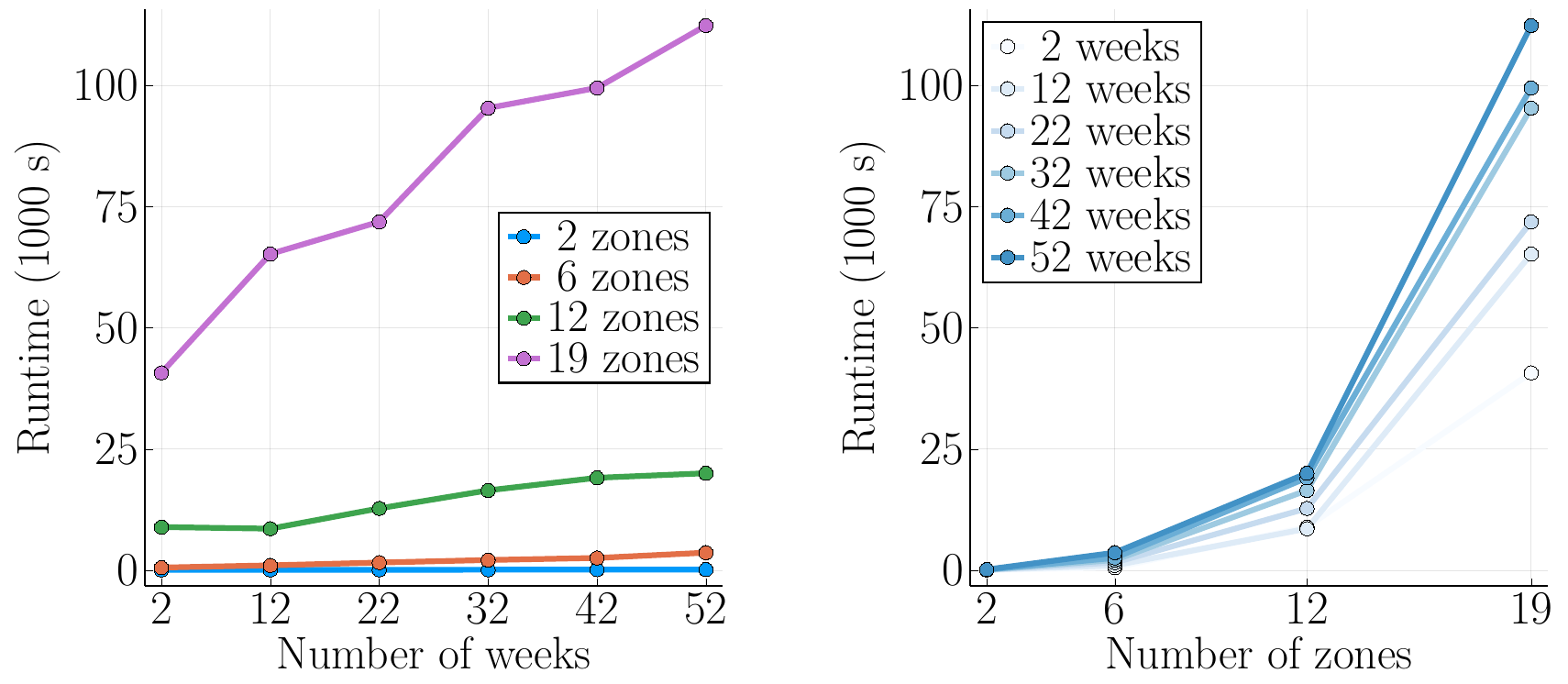}
    \caption{Runtime by weeks (left) and by zone (right) obtained by applying Benders decomposition to solve the MILP problem in the CO$_2$-constrained case. Runtime grows linearly with the number of weeks, while it increases quadratically with the number of zones.}
    \label{fig:runtime_by_zone_weeks}
\end{figure}

Figure~\ref{fig:runtime_by_zone_weeks} suggests that the runtime associated with our Benders algorithm grows linearly with the number of weeks considered. We observe that runtime grows quadratically with the number of modeled zones. This fact hints that there may be value in decomposing along spatial dimensions in future work, as the algorithm described here operates only over temporal dimensions.

\begin{figure}[h!]
\centering
\includegraphics[width=0.8\textwidth]{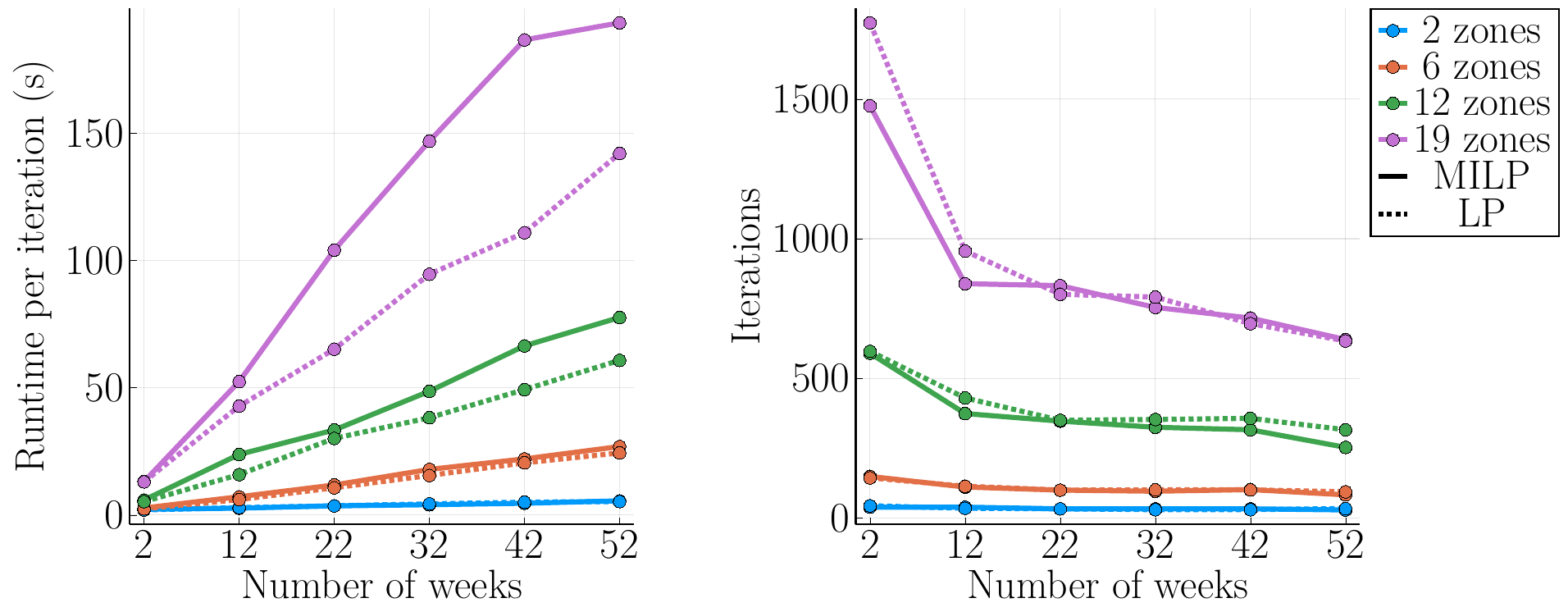}
\caption{Runtime per iteration by week (left) and number of iterations (right), shown for different numbers of zones in the reference case. Runtime per iteration increases both with the number of weeks and the number of zones modeled. Number of iterations increases with the number of zones but decreases with the number of weeks modeled.} \label{runtime_components}
\end{figure}

As the number of weeks in the simulation increases, runtime per iteration increases linearly - see Figure~\ref{runtime_components}. Because the master model incorporates more information per iteration due to its increased number of cuts, the number of iterations required for convergence decreases roughly logarithmically (Figure~\ref{runtime_components}). These trends explain the superior performance of the decomposition method for problems with a larger numbers of weeks.

In contrast, increasing the number of zones increases the size of both the master~\eqref{eq:master} and the subproblem~\eqref{eq:subprob}; accordingly, we see a simultaneous increase in iteration time and number of iterations required for convergence (Figure~\ref{runtime_components}). These dual impacts explain the quadratic runtime increase shown in Figure~\ref{fig:runtime_by_zone_weeks} as the number of zones increases. Increasing the number of zones increases number of resource operational decisions and thus increases the amount of information that must be incorporated by a given cut, as the vector $\bm{y}$ in \eqref{eq:master_cut} grows in dimensionality. The direct relationship between zones and number of iterations suggests that higher-dimensionality cuts poorly approximate the subproblems' objective values.
Our decomposition method outperforms the monolithic approach on all MILP formulations with more than 2 weeks. It outperformed LP formulations on 2- and 6-zone problem instances with more than 12 or 32 weeks, respectively. The decomposition method remained within an order of magnitude of runtime on all simulations and always maintained a linear relationship between weeks and runtime, suggesting that the decomposed model will eventually outperform a monolithic approach as temporal resolution increases.

We further investigate the benefits of our novel decomposition scheme when compared to a standard implementation of Benders decomposition algorithm, which keeps timesteps coupled when solving Problem~\eqref{eq:prob}. In this case, the master problem is given by equation ~\eqref{eq:master_classic}, while the single operational subproblem is formulated as in~\eqref{eq:subprob_classic}.

The Benders algorithm with a fully coupled operational subproblem is intractable under all instances with more than 22 weeks of operation (Table~\ref{fullop_table}). An aggregated Benders cut does not reduce the number of iterations as temporal resolution increases and further slows iteration time due to the size of the subproblem and the impossibility of parallelizing it. 
\begin{table}[!htbp]
\centering
\footnotesize{
 \renewcommand{\arraystretch}{1.45}
    \caption{Runtime in seconds and number of iterations for 6-Zone LP CO$_2$-constrained simulations. 6-Zone CO$_2$ was selected for representation because it was the largest case with tractable full operation models. Intractable simulations are marked {\scshape INT}.}
    \label{fullop_table}
    \begin{tabular}{r|cccccc|cccccc}
    \toprule
    \multicolumn{1}{c}{Model Type} & \multicolumn{6}{c}{Runtime (100s)} & \multicolumn{6}{c}{Iterations}\\
    \cmidrule(lr){1-1}\cmidrule(lr){2-7}\cmidrule(lr){8-13}
    Weeks $\rightarrow$ & 2 & 12 & 22 & 32 & 42 & 52 & 2 & 12 & 22 & 32 & 42 & 52 \\
    \midrule
    Monolithic & 0.6 & 4.3 & 10.8 & 18.3 & 27.5 & 36.9 & \multicolumn{6}{c}{NA} \\
    Benders & 5.9 & 10.1 & 16.6 & 19.7 & 26.0 & 28.3 & 231 & 139 & 132 & 104 & 107 & 94 \\
    Benders Full-Op & 40.6 & 5.7e2 & 1.5e3 & {\scshape INT} & {\scshape INT} & {\scshape INT} & 509 & 607 & 630 & {\scshape INT} & {\scshape INT} & {\scshape INT} \\
    \bottomrule
    \end{tabular}}
\end{table}

\begin{figure}[h!]
\centering
    \includegraphics[width=0.5\textwidth]{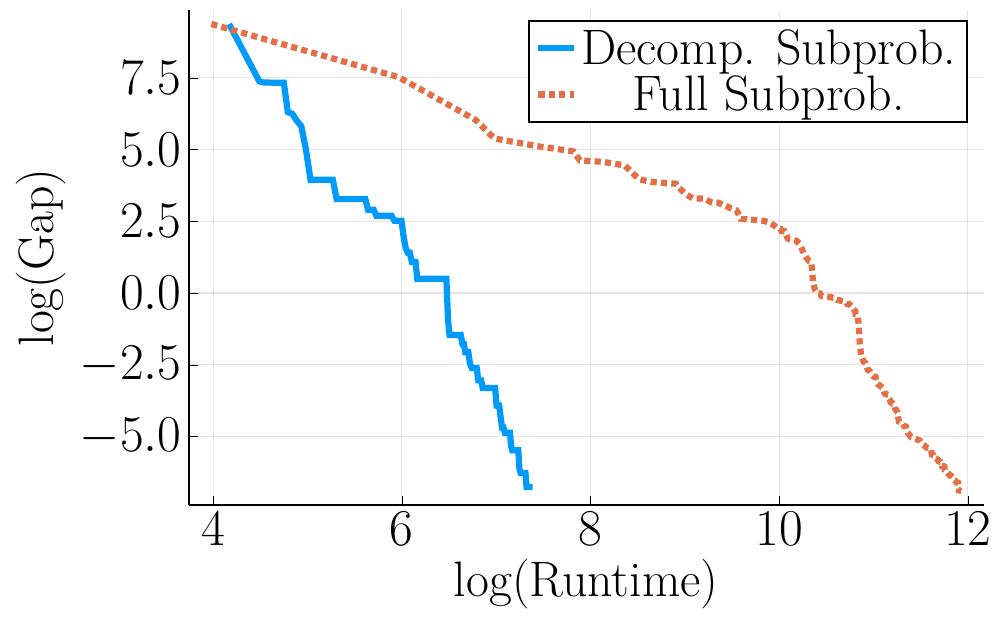}
    \caption{Optimality gap vs. runtime for our Benders algorithm with decomposed subproblem~\eqref{eq:subprob} and a standard Benders implementation with a full operational subproblem~\eqref{eq:subprob_classic}. Plotted on logarithmic x- and y-axes.} 
    \label{fullop_runtime}
\end{figure}

Because it is the largest simulation found to be tractable under this standard Benders implementation, we consider a 6-zone LP network run on 22 weeks of data under the CO$_2$ policy scenario. Figure \ref{fullop_runtime} compares the performance of the decomposition scheme developed in Section~\ref{sec:solution_method} with the standard Benders implementation as described in equations~\ref{eq:master_classic}~and~\ref{eq:subprob_classic}. For both algorithms, we define the optimality gap as $\frac{\text{UB}-\text{LB}}{\text{LB}}$, where $\text{UB}$ and $\text{LB}$ are current upper and lower bounds on the optimal value, respectively. Figure~\ref{fullop_runtime} shows that the novel decomposition algorithm takes significantly less time to converge, running in 132 iterations with an average iteration time of 12.0 seconds. The standard Benders algorithm, by contrast, runs in 630 iterations with an average iteration time of 241.0 seconds (see Table \ref{fullop_table}). 

\label{MSE_explanation}
The purpose of the energy systems planning models is to provide decision support for resource capacity expansion. As such, the specific solutions vector for these capacity planning problems is also of salient interest to stakeholders and decision makers (e.g. in utility integrated resource planning processes.) Because system cost is not the primary output of these models, we do not believe it appropriate to rely on it as the sole indicator of the quality of models' results. In fact, \cite{palmintier2013heterogeneous} note that error in cost does not correlate with error in other system characteristics depending on the means of clustering used or capacity installed for different resources in the system. It is thus impossible to extrapolate error in capacity decisions based solely on error in cost. Therefore, we define the mean squared error (MSE) in investment decisions as:
\begin{equation}
\label{eq:mse}
    \text{MSE}(j) =\frac{1}{|G|}\sqrt{\sum_{g \in G}({\pCAP_g}^j - {\pCAP_g}^{52})^2},
\end{equation}
where ${\pCAP_g}^j$ is the recommended capacity for resource $g \in G$, computed solving our planning problem for $j \in \{2,12,22,32,42,52\}$ weeks. 

These MSE values are computed by comparing near-optimal solutions for the considered MILPs. There may be multiple investment portfolios yielding similar costs for any given model. The MSE metric described in Equation~\eqref{eq:mse} assesses the ability of lower-resolution models to recapitulate investment recommendations of higher-resolution models which inherently have lower structural uncertainty, even considering the possibility of multiple near-optimal solutions.

\begin{figure}[h!]
\centering
\includegraphics[width=0.5\textwidth]{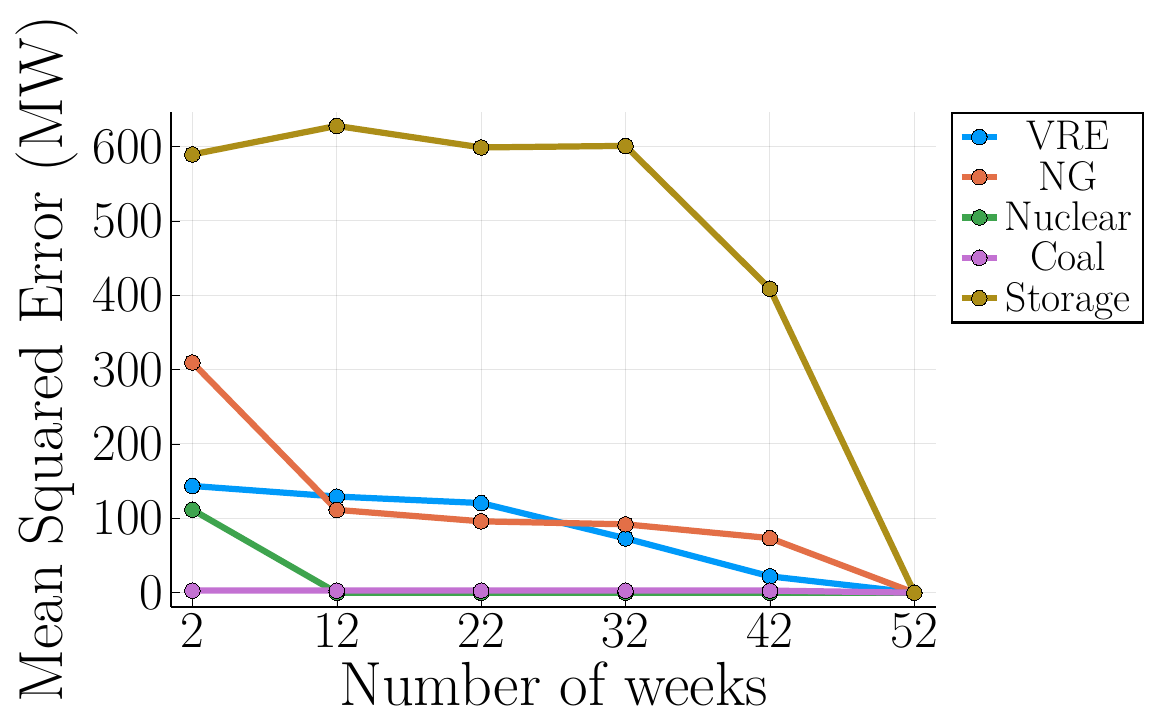}
\caption{Mean Squared Error as defined in Equation~\eqref{eq:mse}. Figure shows that low temporal resolution causes great deviation in capacity recommendations for storage, NG, and VRE. Data is shown for the 19-zone, CO$_2$-constrained MILP trial.}
\label{fig:MSE}
\end{figure}

Figure \ref{fig:MSE} shows that MSE decreases as the temporal resolution increases and that using too few representative weeks can lead to average deviations as high as 600 MW per site for storage clusters. While we show that MSE decreases as temporal resolution increases, it is impossible to compute rigorous upper bounds on MSE \emph{a priori}. That is, it is impossible to prove the maximum deviation from optimal capacity investments due to abstractions without solving the full model for comparison.

Figure~\ref{fig:MSE} shows that temporal resolution has greatest impact on storage, NG, and VRE resources, likely due to the misrepresentation of VRE availability in models with low temporal resolution. The relatively high MSE values for resources in Figure~\ref{fig:MSE} demonstrate models' difficulty in providing recommendations on storage, NG, and VRE investments at high levels of temporal aggregation.

Coal resources are incentivized to retire completely by policy constraints. Similarly, nuclear resources are disincentivized from building due to a high fixed cost of investment. These impacts yield little variation in capacity regardless of temporal resolution, leading to low MSE in Figure~\ref{fig:MSE}. Errors at the resource level are likely to vary with policy or cost assumptions.

\section{Conclusions}
Our algorithm was able to solve MILP energy systems planning problems that were intractable when monolithic; we showed that this ability markedly decreased error associated with modeling. Furthermore, Figure~\ref{fullop_runtime} shows that the standard Benders algorithm without decoupling via budgeting constraints in the master model is significantly less successful than our fully decomposed scheme, often leading to intractibility on smaller problems than the monolithic approach (Table~\ref{fullop_table}.) Standard decomposition schemes that do not separate timesteps do not gain as much information from Benders cuts at each iteration. We therefore note that the reformulation of the problem from \eqref{eq:prob} to \eqref{eq:prob_decomp} is the primary contribution of this paper.

The proposed decomposition scheme has several advantages over monolithic models currently in use industry-wide and methods explored in previous literature where operational subproblems were not decomposed~\citep{lara2018deterministic,li2022mixed} or were decomposed via a nested algorithm~\citep{sepulveda2020decarbonization}. Below, we list some specific benefits to the proposed algorithmic scheme. We note as well that model improvements are not restricted to the framework of energy systems and can inform researchers working on integrated planning problems in other application areas, including water resources~\citep{naderi2017}, industrial processes~\citep{shah2012}, and facility location~\citep{an2020}.
\begin{enumerate}
    \item \textit{Superior Performance:} Our decomposed model consistently outperformed monolithic solution approaches using state-of-the-art commercial solvers on cases with discrete, integer investment decisions and was competitive with linearized investment decisions depending on model size and structure tested. Our ability to decouple and parallelize subperiods and provide $\vert W \vert$ cuts per iteration decreased both the number of iterations and runtime per iteration relative to a more conventional Benders decomposition algorithm, resulting in linear runtime increase as a function of resolution.
    \item \textit{Increased Resolution:} Because runtime scales well with master model size, a great deal of information on investments (including geographic constraints for expansion, integer transmission expansion, and different constraints or properties for different units of thermal plant) can be included with minimal expense to performance. Our decomposition algorithm's runtime scales linearly with number of weeks included, which further allows for inclusion of more subperiods beyond the standard strategy of considering only few weeks or days of operation without risking intractability. Decreasing the level of abstraction in this way helps eliminate structural uncertainty~\citep{pfenninger2017dealing}. Indeed, we demonstrated that increasing temporal resolution decreased the MSE associated with investment decisions by resource type and location.
\end{enumerate}

\subsection{Future Directions and Novel Capabilities}
While increased resolution (as a corollary, decreased reliance on abstraction) is one exciting implication of this work, it is not the only application of the proposed algorithmic scheme. Improved performance also benefits research due to the additional analysis enabled by decreased computational bandwidth. Some potential additional applications include:

\begin{enumerate}
    \item \textit{Enabling More Extensive Analysis:} Improved computational performance due to decomposition can also be employed to conduct more extensive exploration of parameter uncertainty via scenario analysis methods or enable incorporation of methods like multi-objective optimization or modeling to generate alternatives~\citep{patankar2020land}, which explores a wide range of alternative solutions with similar objective function results (e.g. costs). Such methods require solving the planning problem multiple times under different parameters or objective functions, thus benefiting greatly from improved runtime for each problem.
    \item \textit{Ability to Capture Economies of Scale:} Making integer investments tractable permits modeling of discrete capacity investment decisions that capture economics of unit scale that are common in many applications, including electricity transmission and generator investment decisions. Other methods with continuous capacity decisions cannot capture economies of unit scale and may result in  unrealistic biases~\citep{donohoo2014design}. Additionally, without decomposition, modeling electricity network and power flows with KVL introduces nonlinearities in combined investment / operations models and generally requires Big-M reformulations that can significantly slow computational time loosening the convex relaxation of the MILP problem. The structure of the decomposition algorithm herein separates investment and operations into discrete problems and thus can allow KVL and other operational characteristics that have interactions with investments to be included seamlessly in a convex model.
    \item \textit{Model Accessibility.} Current state-of-the-art energy systems planning tools are inaccessible to many potential users due to high computational demands and their associated need for expensive commercial solvers. In fact, the benchmarks computed by Hans Mittelmann~\citep{mittelmann2023} report that the fastest open source solver can still be, on average, 30 to 40 times slower than the best commercial solvers. Analogously, \cite{han2021comprehensive} show that solvers that are completely open source (score ``10" on the ``openness" scale) tend to perform poorly on large-scale security constrained economic dispatch models. The decomposition method introduced here involves solving substantially smaller master and operational subproblems, which makes it easier to implement open-source solvers, thereby increasing accessibility of macro-energy systems planning models.
\end{enumerate}

\subsection*{Acknowledgments}
Funding for this work was provided by the Princeton Carbon Mitigation Initiative (funded by a gift from BP) and the Princeton Zero-carbon Technology Consortium (funded by gifts from GE, Google, ClearPath, and Breakthrough Energy).
\newpage
\appendix
\section{Notation}
\label{sec:notation}
\begin{table}[h!]
\centering
\footnotesize{
\renewcommand{\arraystretch}{1.45}
\caption{A list of all sets included in the simulation.} \label{nomenclature_sets}
\begin{tabular}{p{0.08\textwidth}p{0.36\textwidth}|p{0.08\textwidth}p{0.36\textwidth}}
\toprule
\multicolumn{1}{c}{\scshape Set} & \multicolumn{1}{c}{\scshape Composes} & \multicolumn{1}{c}{\scshape Set} & \multicolumn{1}{c}{\scshape Composes} \\
\midrule
$G$ & All resources & $\STOR$ & All storage resources \\
$\UC$ & All resources subject to unit commitment & $\HYDRO$ & All hydro power resources \\
$G^{NONRET}$ & All resources that cannot be retired & $S$ & All consumer segments (demand)\\
$\RPS$ & All resources qualifying for RPS policy & $L$ & All transmission lines \\
$Z$ & All spatial zones  & $W$ & All subperiods \\
$H_w$ & All hours per subperiod $w \in W$ & $\OutboundLines_z$ & Lines carrying power out of zone $z$ \\
$\InboundLines_z$ & Lines carrying power into zone $z$ & & \\
\bottomrule
\end{tabular}
}
\end{table}
\begin{table}[h!]
\centering
\footnotesize{
\caption{A list of all parameters in the simulation.}
 \renewcommand{\arraystretch}{1.75}
\begin{tabular}{m{0.05\textwidth}m{0.4\textwidth} | m{0.05\textwidth}m{0.4\textwidth}}
\toprule
\multicolumn{1}{c}{\scshape Param} & \multicolumn{1}{c}{\scshape Definition} & \multicolumn{1}{c}{\scshape Param} & \multicolumn{1}{c}{\scshape Definition}
\\\cmidrule(lr){1-2}\cmidrule(lr){3-4}
\multicolumn{4}{c}{\scshape Investments and Capacity} \\
\midrule
 $\pCapEx_g$ & Existing capacity [MW], resource $g$ & $\pCapMax_g$ & Max capacity [MW], resource $g$\\
$\eCapEx_g$ & Existing storage capacity [MWh], resource $g$ &  $\eCapMax_g$ & Max storage capacity [MWh], resource $g$ \\
$\pCapSize_g$ & Capacity size [MW]$^\dagger$, resource $g$ & $\eCapSize_g$ & Storage capacity size [MWh]$^\dagger$, resource $g$\\
 $\tCapEx_l$ & Existing transmission capacity [MW], line $l$ & $\tMaxTransPossible_l$ & Max transmission capacity [MW], line $l$ \\
$\MinDuration_g$ & Min duration, resource $g^\ast$ [MWh/MW] & $\MaxDuration_g$ & Max duration, resource $g^\ast$ [MWh/MW] \\

\midrule
\multicolumn{4}{c}{\scshape Operations} \\
\midrule

$\variability_{g,t}$ & Availability [\%], resource $g$, timestep $t$ & $d_{z,t}$ & Net demand [MWh], zone $z$, timestep $t$\\
$\NseMax_{s}$ & Max non-served energy (NSE) [\%], segment $s$ & $\MinPower_g$ & Min output [\%], resource $g$\\
$\periodweight_w$ & Weight assigned to subperiod $w$ [-] & $\hoursperperiod_w$ & Number of hours in subperiod $w$ [-]\\
$\sampleweight_t$ & Weight assigned to hour $t$ [-]&  $\Duration_g$ & Duration$^\ast$ for hydro resource $g$ [MWh/MW]\\
$\EffUp_g$ & Charging efficiency [\%], storage resource $g$ & $\EffDown_g$ & Discharging efficiency [\%], storage resource $g$  \\
$\SelfDisch_g$ & Self-discharge rate from storage resource $g$ [\%]& $\NormHourlyInflow_{g,t}$ & Norm. inflow, hydro resource $g$, timestep $t$ [\%]  \\
$\RampUp_g$ & Max ramp up rate [\%/hr], resource $g$ & $\RampDown_g$ & Max ramp down rate [\%/hr], resource $g$ \\
$\UpTime_g$ & Min up time [hours], resource $g$ & $\DnTime_g$ & Min down time [hours], resource $g$\\
\midrule
\multicolumn{4}{c}{\scshape Policies}\\
\midrule
$\rpscap$ & Share of demand in RPS constraint [\%] & $\carboncap$ & CO2 emission cap [tons/MWh] \\
$\carbonemiss_g$ & CO2 emission factor [tons/MWh], resource $g$ & \\
\midrule
\multicolumn{4}{c}{\scshape Costs}\\
\midrule
$\InvCostMWyr_g$ & Cost of investment in resource $g$ [\$/MW-yr] & $\FOMCostMWyr_g$ & Fixed O\&M cost of resource $g$ [\$/MW-yr] \\
$\InvCostMWhyr_g$ & Cost of investment, storage resource $g$ [\$/MWh-yr] & $\FOMCostMWhyr_g$ & Fixed O\&M cost, storage resource $g$ [\$/MWh-yr] \\
$\LineCostMWyr_l$ & Cost of investment in line $l$ [\$/MW-yr] & $\VarCost_{g,t}$ & Variable costs [\$/MWh], resource $g$, timestep $t$\\
$\NSECost_{s,z}$ & Cost of NSE [\$/MWh], segment $s$, zone $z$ & $\StartUpCost_{g}$ & Cost to start up resource $g \in \UC$ [\$] \\
$\RPSCost$ & Cost of RPS constraint noncompliance [\$/MWh] & $\CarbonCapCost$ & Cost of CO$_2$ constraint noncompliance [\$/tons] \\
\bottomrule
\multicolumn{4}{c}{\footnotesize $\ast$ ``Duration'' here refers to the ratio between power and energy of a given resource, in MWh/MW}
\\ \multicolumn{4}{c}{\footnotesize $\dagger$ Set to 1 for $g \notin \UC$}
\end{tabular}
}
\end{table}
\begin{table}[h!]
\centering
\footnotesize{
\caption{A list of all variables in the simulation.}
 \renewcommand{\arraystretch}{1.45}
\begin{tabular}{m{0.05\textwidth}m{0.4\textwidth} | m{0.05\textwidth}m{0.4\textwidth}}
\renewcommand{\arraystretch}{1.45}\\
\toprule
\multicolumn{1}{c}{\scshape Var.} & \multicolumn{1}{c}{\scshape Definition} & \multicolumn{1}{c}{\scshape Var.} & \multicolumn{1}{c}{\scshape Definition}
\\\cmidrule(lr){1-2}\cmidrule(lr){3-4}
\multicolumn{4}{c}{\scshape Investments and Capacity} \\
\midrule
$\pCAP_g$ & Capacity [MW], resource $g$ & $\pNEWCAP_g$ & Investments in generation resource $g^{*}$ [-] \\
$\eCAP_g$ & Capacity [MWh], storage resource $g$ & $\eNEWCAP_g$ & Investments in storage resource $g^{*}$ [-] \\
$\tCAP_l$ & Capacity$^*$ [MW], transmission line $l$ & $\tNEWCAP_l$ & Investments in transmission line $l^{*}$ [MW]   \\
$\pRETCAP_g$ & Retirements, generation resource $g^{*}$ [-]& $\eRETCAP_g$ & Retirements, storage resource $g^{*}$ [-]  \\
\midrule
\multicolumn{4}{c}{\scshape Operations} \\
\midrule
$\pGEN_{g,t}$ & Generation[MWh], resource $g$, timestep $t$ & $\pCHARGE_{g.t}$ & Withdrawals [MWh], of $g \in \STOR$, timestep $t$ \\
$\pNSE_{s,z,t}$ & NSE [MWh], segment $s$, zone $z$, timestep $t$ & $\eSOC_{g,t}$ & State of charge [MWh] for $g \in \STOR$, timestep $t$\\
$\eLEVEL_{g,t}$ & Reservoir level [MWh] for $g \in \HYDRO$, timestep $t$ & $\pSPILL_{g,t}$ & Spillage [MWh] from $g \in \HYDRO$ in timestep $t$ \\
$\pFLOW_{l,t}$ & Flow [MWh] across line $l$ in timestep $t$ & $\vCOMMIT_{g,t}$ & Units [-] of $g \in \UC$ commited in timestep $t$ \\
$\vSHUT_{g,t}$ & Units [-] of $g \in \UC$ shut down in timestep $t$ [-] & $\vSTART_{g,t}$ & Units [-] of $g \in \UC$ started up in timestep $t$ \\
\midrule
\multicolumn{4}{c}{\scshape Policies} \\
\midrule
$\vRPS_w$ & Noncompliance [MWh] with RPS policy & $\vCarbonCap_w$ & Noncompliance [tons / MWh] with CO2 cap policy\\
\bottomrule
\multicolumn{4}{c}{\footnotesize $\ast$ Variable is included as integer in the MILP problem formulations}
\end{tabular}}
\end{table}

\cleardoublepage
\section{Problem Description} \label{prob_description}
To complement the constraints shown in the main text as a monolithic model, we present below an in-depth description of the separate groups of constraints in our model formulation.

For descriptions of variables, sets, and parameters used in sections \ref{prob_description_inv}, \ref{prob_description_op}, \ref{prob_description_policy}, and \ref{prob_objective}, refer to the nomenclature table in section~\ref{sec:notation}.

\subsection{Investment constraints} \label{prob_description_inv}

The following constraints model resource capacity given investments and retirements:
\begin{subequations}
\label{eq:invest_cons_1}
\begin{align}
& \pCapSize_{g}\pNEWCAP_g \leq \pCapMax_{g}, \quad g \in G \\
& \pCapSize_{g}\pRETCAP_g \leq \pCapEx_{g}, \quad g \in G \\
& \pRETCAP_g = 0, \quad g \in G^{NONRET} \\
& \pCAP_g = \pCapEx_g +\pCapSize_{g} (\pNEWCAP_g - \pRETCAP_g), \quad g \in G,
\end{align}
\end{subequations}

Additional variables and constraints are included to allow for the modeling of energy storage resources. The following constraints model investments in storage:
\begin{subequations}
\label{eq:invest_cons_2}
\begin{align}
& \eCapSize_{g}\eNEWCAP_g \leq \eCapMax_{g}, \quad  g \in \STOR \\
& \eCapSize_{g}\eRETCAP_g \leq \eCapEx_{g}, \quad  g \in \STOR \\
& \eRETCAP_g = 0, \quad  g \in G^{NONRET} \cap {\STOR}\\
& \eCAP_g = \eCapEx_g + \eCapSize_{g}(\eNEWCAP_g - \eRETCAP_g), \quad  g \in \STOR\\
& \MinDuration_g\pCAP_g  \leq \eCAP_g, \quad  g \in \STOR\\
& \MaxDuration_g\pCAP_g  \geq \eCAP_g, \quad  g \in \STOR,
\end{align}
\end{subequations}

Transmission capacity expansion is modeled by:
\begin{subequations}
    \label{eq:invest_cons_3}
    \begin{align}
        & \tNEWCAP_l \leq \tMaxTransPossible_l, \quad  l \in L \\
        & \tCAP_l = \tCapEx_l + \tNEWCAP_l, \quad  l \in L,
    \end{align}
\end{subequations}

\subsection{Operational constraints}
\label{prob_description_op}

System operation is primarily subject to the power demand balance stating that zonal demand equals available power within a zone (the sum of generation and net storage discharge and transmission imports) plus curtailed demand:
\begin{equation}
\label{eq:dem_balance}
\sum_{g \in G_z} \pGEN_{g,t} - \sum_{g \in \STOR_z} \pCHARGE_{g,t}  -\sum_{l \in \OutboundLines_z}\pFLOW_{l,t} +\sum_{l \in \InboundLines_z}\pFLOW_{l,t}  + \sum_{s\in S}\pNSE_{s,z,t} = d_{z,t},\quad  z\in Z, \; t\in H_w, \;  w \in W,
\end{equation}

Some resources have a constant available dispatchable capacity equal to their installed capacity, while other resources (e.g., VRE) may vary in maximum output throughout the planning period (e.g. based on wind or solar resource variability). To capture this, we introduce a maximum power output parameter $\sigma_{g,t} \in [0,1]$ to model the fraction of installed power capacity available for dispatch at time $t \in H_{w}$. Note that for many resources (e.g., thermal and storage,) $\sigma_{g,t} = 1$ to denote the resource being consistently available at its full capacity. We also consider parameters $\Duration_g$ representing the ratio of power to energy for storage in hydropower resource $g \in \HYDRO$, $\EffUp_g$ representing the charging efficiency for storage resource $g \in \STOR$, and $\EffDown_g$ representing the discharging efficiency for storage resource $g \in \STOR$. Dispatched power and stored energy are limited by installed capacity, resulting in the following linear constraints:
\begin{subequations}
\label{eq:operation_cons_1}
\begin{align}
    &\pGEN_{g,t} \leq \variability_{g,t}\pCAP_g , \quad  g\in G \setminus \UC, \;  t \in H_w, \; w \in W\\
    &\pCHARGE_{g,t} \leq \variability_{g,t}\pCAP_g , \quad  g \in \STOR,  \;  t \in H_w, \; w \in W\\
    &\pGEN_{g,t} + \pCHARGE_{g,t} \leq \pCAP_g, \quad  g \in \STOR, \;  t \in H_w, \; w \in W\\
    &\eSOC_{g,t} \leq \eCAP_g, \quad  g \in \STOR, \;  t \in H_w, \; w \in W\\
    &\eLEVEL_{g,t} \leq \Duration_g\pCAP_g, \quad  g \in \HYDRO, \;  t \in H_w, \; w \in W\\
    &\EffUp_g \pCHARGE_{g,t} \leq \eCAP_g - \eSOC_{g,t}, \quad  g \in \STOR, \;  t \in H_w, \; w \in W\\
    & \frac{\pGEN_{g,t}}{\EffDown_g} \leq \eSOC_{g,t}, \quad  g \in \STOR, \;  t \in H_w, \; w \in W.
\end{align}
\end{subequations}
Some resources have a minimum required output level, which is enforced by:
\begin{subequations}
\label{eq:operation_cons_2}
\begin{align}
    & \pGEN_{g,t} \geq \MinPower_g \pCAP_g, \quad  g \in G \setminus(\UC \cup \STOR \cup \HYDRO), \;  t \in H_w, \; w \in W\\
    & \pGEN_{g,t} + \pSPILL_{g,t} \geq \MinPower_g \pCAP_g, \quad  g \in \HYDRO, \;  t \in H_w, \; w \in W.
\end{align}
\end{subequations}
Power flow between zones is bounded by transmission capacity:
\begin{subequations}
\label{eq:operation_cons_3}
\begin{align}
    & \pFLOW_{l,t} \leq \tCAP_l, \quad  l\in L, \;  t \in H_w, \;  w \in W\\
    & -\pFLOW_{l,t} \leq \tCAP_l, \quad  l\in L, \;  t \in H_w, \;  w \in W.
\end{align}
\end{subequations}
The maximum curtailed demand in each consumer segment is constrained to be at most a fraction $\NseMax_s \in [0,1]$ of the zonal demand:
\begin{equation}
\label{eq:operation_cons_4}
\pNSE_{s,z,t} \leq \NseMax_s d_{z,t}, \quad s \in S, z\in Z, \;  t \in H_w, \;  w \in W.
\end{equation}
Next, we introduce notation to model storage, UC and ramping limits within subperiods. The first time step in a subperiod is denoted by $t^0_w = (w-1)\hoursperperiod_w+1$, for all $w \in W$, while the last time step in a subperiod is $t_w = \hoursperperiod_w w$, for all $w \in W$. The time steps that are not at the start of a subperiod are included in subset $H^0_w = H_w \setminus \{t^0_w\}$, for all $w \in W$.
Storage and hydropower operational constraints are as follows:
\begin{subequations}
\label{eq:operation_cons_5}
\begin{align} 
& \eSOC_{g,t} -\eSOC_{g,t-1} = \EffUp_g \pCHARGE_{g,t} - \frac{\pGEN_{g,t}}{\EffDown_g} - \SelfDisch_g \eSOC_{g,t-1},\quad g \in \STOR, \; t \in H_w^0, \; w \in W\\
& \eSOC_{g,t^0_w} -\eSOC_{g,t_w} = \EffUp_g \pCHARGE_{g,t^0_w} - \frac{\pGEN_{g,t^0_w}}{\EffDown_g} - \SelfDisch_g \eSOC_{g,t_w},\quad g \in \STOR, \; w \in W\\
& \eLEVEL_{g,t} -\eLEVEL_{g,t-1} =  \NormHourlyInflow_{g,t} \pCAP_g - \pGEN_{g,t} - \pSPILL_{g,t},\quad g \in \HYDRO, \; t \in H_w^0,\; w \in W\\
& \eLEVEL_{g,t^0_w} -\eLEVEL_{g,t_w} =  \NormHourlyInflow_{g,t^0_w} \pCAP_g - \pGEN_{g,t^0_w} - \pSPILL_{g,t^0_w},\quad g \in \HYDRO, \; w \in W.
\end{align}
\end{subequations}
Ramping limits for thermal units are enforced within each subperiod. The difference in discharge from a resource $g$ in two consecutive steps is constrained as follows:
\begin{subequations}
\label{eq:operation_cons_6}
\begin{align}
& \pGEN_{g,t} - \pGEN_{g,t-1}  \leq \RampUp_g \pCAP_g, \quad g \in G \setminus \UC, \; t\in H_w^0,\; w \in W \\
& \pGEN_{g,t-1} - \pGEN_{g,t}  \leq \RampDown_g \pCAP_g, \quad g \in G \setminus \UC, \; t\in H_w^0,\; w \in W  \\
& \pGEN_{g,t^0_w} - \pGEN_{g,t_w}  \leq \RampUp_g \pCAP_g, \quad g \in G \setminus \UC,\; w \in W\label{eq:operation_cons_6_ramp_up} \\
& \pGEN_{g,t_w} - \pGEN_{g,t^0_w}  \leq \RampDown_g \pCAP_g, \quad g \in G \setminus \UC,\; w \in W .\label{eq:operation_cons_6_ramp_down}
\end{align}
\end{subequations}
Resources in $\UC$ can commit an integer number of units at each time step. To avoid having integer operational variables, we relax integrality constraints on UC variables. This is a common strategy in the solution of energy systems planning problems \citep{lara2018deterministic,li2022mixed}, as integer variables are expected to have tight linear relaxations for UC operations - see the optimality gaps reported in Table 2 in \cite{li2022mixed}. Note that we do not relax integrality constraints on investment decision variables, whose integer representation is critical to fully capture economies of scale in generation and transmission expansion. UC is modeled as follows:
\begin{subequations}
\label{eq:operation_cons_7}
\begin{align}
& \pCapSize_g\vCOMMIT_{g,t} \leq \pCAP_g, \quad g \in \UC, \; t \in H_w, \; w \in W \\
& \pCapSize_g\vSTART_{g,t} \leq \pCAP_g, \quad g \in \UC , \; t \in H_w, \; w \in W \\
& \pCapSize_g\vSHUT_{g,t} \leq \pCAP_g, \quad g \in \UC , \; t \in H_w, \; w \in W \\
& \pGEN_{g,t} \geq \vCOMMIT_{g,t} \MinPower_g \pCapSize_g,\quad g \in \UC , \; t \in H_w, \; w \in W \\
& \pGEN_{g,t} \leq \vCOMMIT_{g,t} \variability_{g,t} \pCapSize_g,\quad g \in \UC , \; t \in H_w, \; w \in W \\
& \vCOMMIT_{g,t} - \vCOMMIT_{g,t-1} = \vSTART_{g,t} - \vSHUT_{g,t}, \quad g \in \UC, \; t\in H_w^0, \; w \in W \\
& \vCOMMIT_{g,t^0_w} - \vCOMMIT_{g,t_w} = \vSTART_{g,t^0_w} - \vSHUT_{g,t^0_w}, \quad g \in \UC, \;  w \in W .
\end{align}
\end{subequations}
Ramping limits on resources in $g \in \UC$ are given by:
\begin{subequations}
\label{eq:operation_cons_8}
\begin{align}
&\pGEN_{g,t} - \pGEN_{g,t-1} \leq \begin{aligned}[t]\pCapSize_g\RampUp_g( \vCOMMIT_{g,t} - \vSTART_{g,t}) + \pCapSize_g\min(\variability_{g,t},\max(\MinPower_g,\RampUp_g)) \vSTART_{g,t}\\ - \pCapSize_g\MinPower_g \vSHUT_{g,t} , \quad g \in \UC, \; t\in H_w^0, \; w \in W \end{aligned}\\
&\pGEN_{g,t-1} - \pGEN_{g,t} \leq \begin{aligned}[t] \pCapSize_g\RampDown_g (\vCOMMIT_{g,t} - \vSTART_{g,t}) +\pCapSize_g\min(\variability_{g,t},\max(\MinPower_g, \RampDown_g)) \vSHUT_{g,t} \\ - \pCapSize_g\MinPower_g  \vSTART_{g,t} , \quad g \in \UC, \; h\in H_w^0, \; w \in W \end{aligned}\\
&\pGEN_{g,t^0_w} - \pGEN_{g,t_w} \leq\begin{aligned}[t]  \pCapSize_g\RampUp_g( \vCOMMIT_{g,t^0_w} - \vSTART_{g,t^0_w}) + \pCapSize_g\min(\variability_{g,t^0_w},\max(\MinPower_g,\RampUp_g)) \vSTART_{g,t^0_w} \\ - \pCapSize_g\MinPower_g \vSHUT_{g,t^0_w} , \quad g \in \UC,  w \in W  \end{aligned}\\
&\pGEN_{g,t_w} - \pGEN_{g,t^0_w} \leq \begin{aligned}[t]  \pCapSize_g\RampDown_g (\vCOMMIT_{g,t^0_w} - \vSTART_{g,t^0_w}) + \pCapSize_g\min(\variability_{g,t^0_w},\max(\MinPower_g, \RampDown_g)) \vSHUT_{g,t^0_w} \\ - \pCapSize_g\MinPower_g  \vSTART_{g,t^0_w} , \quad g \in \UC,  w \in W. \end{aligned}
\end{align}
\end{subequations}
In addition, committed resources may be required to stay online (resp. offline) for a minimum period of time $\UpTime$ (resp. $\DnTime)$ before being shut down (resp. restarted). These constraints are modeled within each subperiod via definition of index sets $\CountUpTimes$ and $\CountDnTimes$ in~\eqref{eq:time_definition}, which are then used to incorporate constraints~\eqref{eq:operation_cons_9}, stating:
\begin{subequations}
\label{eq:operation_cons_9}
\begin{align}
&\vCOMMIT_{g,t} \geq \sum_{k \in \CountUpTimes_{g,w}(t)}\vSTART_{g,k}\\
& \vCOMMIT_{g,t} +\sum_{k \in \CountDnTimes_{g,w}(t)}\vSHUT_{g,k} \leq \frac{\pCAP_g}{\pCapSize_g},
\end{align}
\end{subequations}
with
\begin{equation}
\begin{split}
\label{eq:time_definition}
\CountUpTimes_{g,w}(t) &= \{\phi_{w,\UpTime_g}(t),\phi_{w,\UpTime_g-1}(t),\ldots,\phi_{w,1}(t),g_{w,0}(t)=t\}, \quad t \in H_w\\
\CountDnTimes_{g,w}(t) &= \{\phi_{w,\DnTime_g}(t),\phi_{w,\DnTime_g-1}(t),\ldots,\phi_{w,1}(t),g_{w,0}(t)=t\}, \quad t \in H_w\\
\end{split}
\end{equation}
where $\phi_{w,n}(t)$ corresponds to the time index that is $n$ steps before $t$ in $H_w$, where $H_w$ is considered as a circular array of length $\hoursperperiod_w$. For example, if $\hoursperperiod_w=10 \; \forall w$, we have that $H_2=\{11,\ldots,20\}$, and $\phi_{2,1}(11)=20$.
\subsection{Policy constraints} \label{prob_description_policy}
A renewable portfolio standard enforces that a given share, $\rpscap$, of the total demand is produced by qualifying resources, $\RPS \subset G$. In order to avoid infeasiblities, we enforce RPS by assigning a heavy penalty to the violation of the minimum energy share constraint. Let $\vRPS_w \geq 0$ be a slack variable used to evaluate the level of non-compliance with RPS in subperiod $w \in W$. For RPS cases, we include the following constraint:
\begin{subequations}
\label{eq:rps_cons}
\begin{align}
&\sum_{w \in W}\Bigg(\Big(\sum_{t \in H_w}\sum_{g \in \RPS}\sampleweight_t \pGEN_{g,t} \Big) + \vRPS_w\Bigg) \geq \rpscap \Big(\sum_{w \in W} \sum_{t \in H_w}\sum_{z \in Z}\alpha_t d_{z,t}\Big),
\end{align}
\end{subequations}
where $\sampleweight_t$ is the weight assigned to the time step $t$ such that $\sampleweight_t = \frac{\periodweight_w}{\hoursperperiod_w}$, for all $t \in H_w$ and $ w \in W$,  $\periodweight_w$ is the weight associated to representative subperiod $w$ by the timeseries clustering algorithm (i.e., the number of hours represented by subperiod $w$). Observe that $\sampleweight_t=1 \ \forall \ t \in T$ for the full year when $|W|=52$.

As a policy constraint in the CO2 scenario, consider a CO2 emissions cap. Analogously to the RPS constraint, we introduce slack variables $\vCarbonCap_w \geq 0$ to evaluate the level of non-compliance with the CO2 emissions cap in subperiod $w \in W$. In the formulation of scenario CO2, we include the following constraint:
\begin{subequations}
\label{eq:co2_cons}
\begin{align}
&\sum_{w \in W}\Bigg(\sum_{t \in H_w}\Big(\sum_{g \in G}\sampleweight_t \carbonemiss_g\pGEN_{g,t} + \sum_{g \in \STOR}\sampleweight_t \carbonemiss_g \pCHARGE_{g,t}\Big) - \vCarbonCap_w  \Bigg)\leq \carboncap,
\end{align}
\end{subequations}
where $\carbonemiss_g$ represents tons of CO2 emissions per-MWh for resources $g \in G$, and  $\carboncap$ is the maximum emission threshold. In our case, we define:
\begin{equation}
    \carboncap = 0.05 \Big(\sum_{w \in W} \sum_{t \in H_w}\sum_{z \in Z}\alpha_t d_{z,t}\Big).
\end{equation}
For scenario REF, CO2 emissions and generation portfolio are unrestricted. We effectively see $\vRPS_w = 0, \ \vCarbonCap_w = 0, \ \rpscap = 0, \ \carboncap = \infty \quad \forall w \in W$.

\subsection{Objective function} \label{prob_objective}
Planning problems in energy systems minimize over costs of both investment and operations. For our purposes, we annualize investment costs and minimize one net cost function for the operations and investments over the year-long planning period.

Total fixed cost is given by:
\begin{equation}
\label{eq:fixed_cost}
\begin{split}
&\totalFixedCost = \sum_{g \in G} \InvCostMWyr_g \pCapSize_g \pNEWCAP_g  + \sum_{g \in \STOR} \InvCostMWhyr_g \eCapSize_g \eNEWCAP_g + \sum_{g \in \HYDRO} \InvCostMWhyr_g \Duration_g \pNEWCAP_g  \\
&+ \sum_{g \in G} \FOMCostMWyr_g \pCAP_g + \sum_{g \in \STOR} \FOMCostMWhyr_g \eCAP_g + \sum_{g \in \HYDRO} \FOMCostMWhyr_g \Duration_g \pCAP_g + \sum_{l \in L} \LineCostMWyr_l \tNEWCAP_l.
\end{split}
\end{equation}
Next, denote by $\VarCost_g$ resources' variable cost per MWh and compute total variable costs as:
\begin{equation}
\label{eq:var_cost}
\begin{split}
& \VarCost = \sum_{w \in W} \sum_{t \in H_w}\sum_{g \in G} \VarCost_g\sampleweight_t  \pGEN_{g,t} +  \sum_{w \in W} \sum_{t \in H_w}\sum_{g \in \STOR} \VarCost_g\sampleweight_t  \pCHARGE_{g,t}  
\end{split}
\end{equation}
Let $\NSECost_{s,z,t}$ be the cost of curtailing demand in consumer segment $s \in S$, zone $z \in Z$, and time $t \in H_w$ in subperiod $w \in W$. The total cost for curtailed demand, is:
\begin{equation}
\label{eq:nse_cost}
\NSECost = \sum_{w \in W} \sum_{t \in H_w} \sum_{z \in Z} \sum_{s \in S}\NSECost_{s,z}\sampleweight_t\pNSE_{s,z,t}
\end{equation}
Let $\StartUpCost_g$ be cost to start-up a unit of resource $g \in \UC$. The total start-up costs are given by:
\begin{equation}
\label{eq:startup_cost}
\begin{split}
\StartUpCost=\sum_{w \in W} \sum_{t \in H_w}\sum_{g \in \UC}\StartUpCost_g \sampleweight_t \vSTART_{g,t}
\end{split}
\end{equation}
Finally, we include penalty costs for violating policy constraints:
\begin{equation}
\label{eq:policy_cost}
\RPSCost + \CarbonCapCost = \sum_{w \in W}\RPSCost \vRPS_w + \sum_{w \in W}\CarbonCapCost \vCarbonCap_w,
\end{equation}
where $\RPSCost$ is equal to the cost for violating the RPS policy, and $\CarbonCapCost$ is the cost for violating the CO2 emissions cap. Total cost of the problem will be equal to:
\begin{equation}
\label{eq:total_cost}
\totalFixedCost + \VarCost + \NSECost + \StartUpCost + \RPSCost + \CarbonCapCost.
\end{equation}
\cleardoublepage
\subsection{Input Data} \label{powergenome_description}
Data for PowerGenome comes from the Public Utility Data Liberation (PUDL) project and the Annual Technology Baseline (ATB) from the National Renewable Energy Lab (NREL)~\citep{schivley2021powergenome}. Demand and initial capacity for resources is identical across all CO2, RPS, and reference cases. Initial capacity and demand is also identical within each zone for varying spatial extent.
\begin{figure}[h!]
\centering
\caption{Load for the 6-zone, 12-week case. Shows load across all timesteps (\ref{load_all_timesteps}), load for the first week of the simulation (\ref{load_first_timestep}), and a map of the zones included (\ref{6z_map}.)}
\subfloat[Load for the 6-zone 12-week case, first timestep]
{\includegraphics[width=0.4\textwidth]{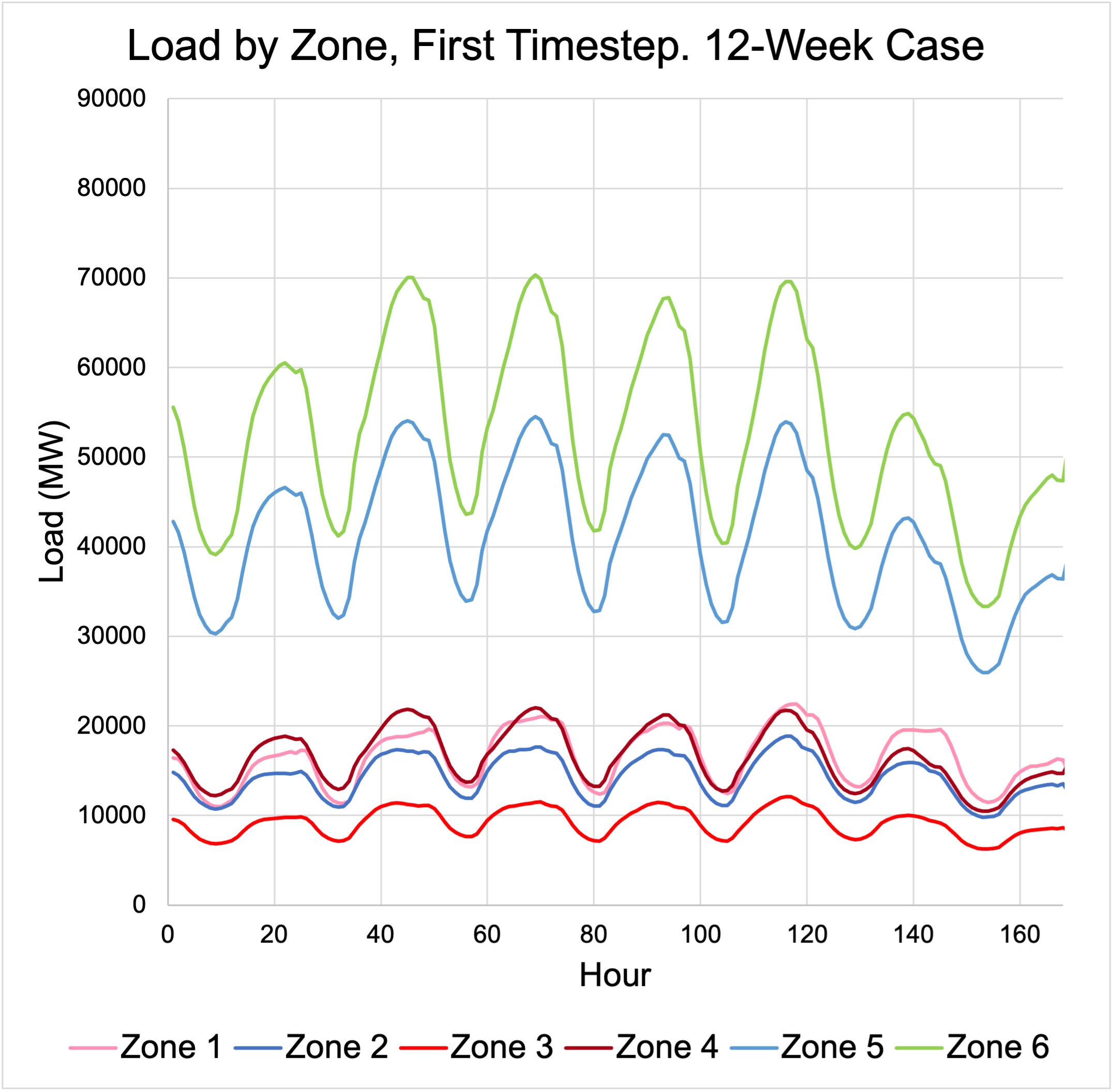}\label{load_first_timestep}}
\qquad
\subfloat[Map of zones in the 6-zone case]{\includegraphics[width=0.3\textwidth]{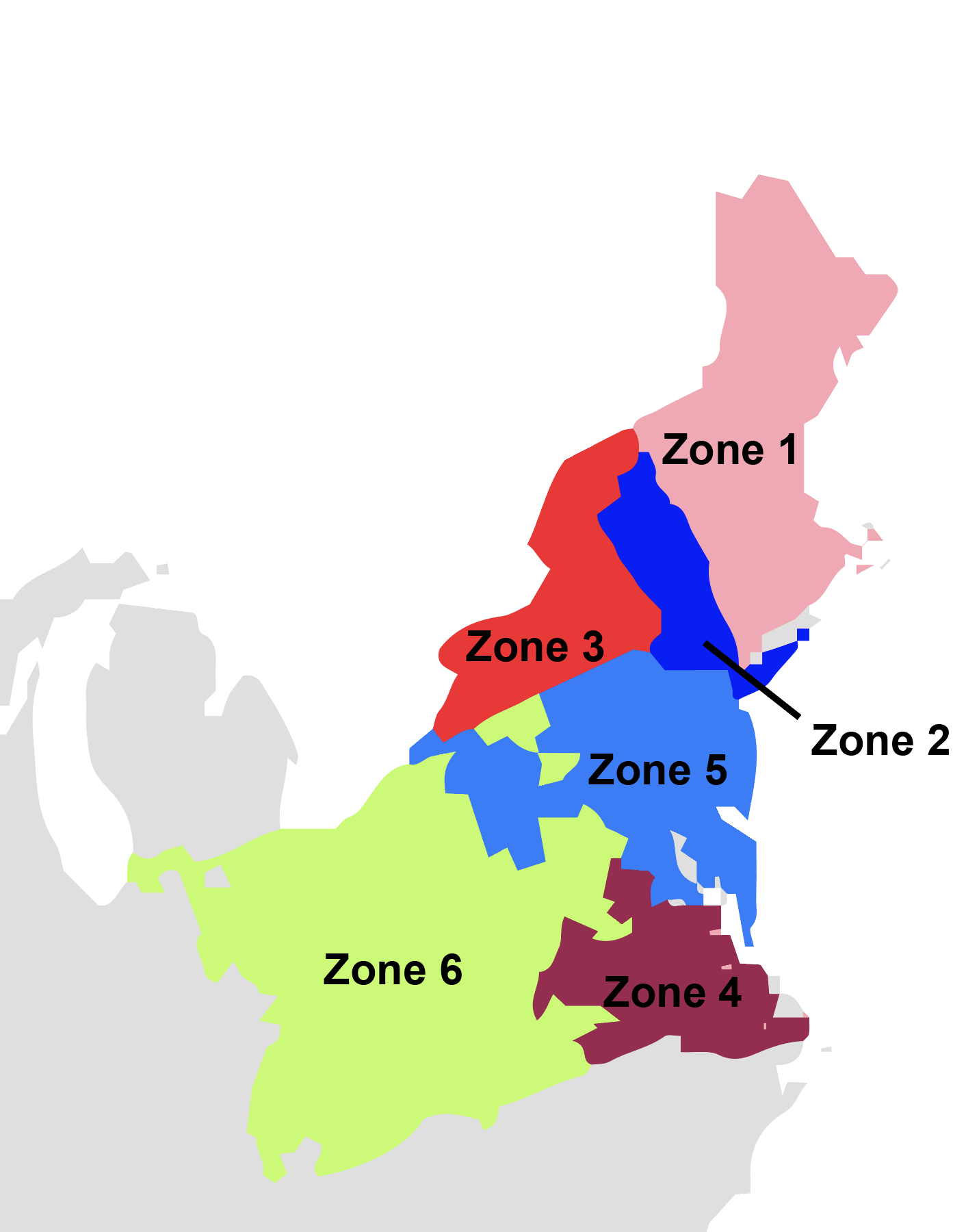}\label{6z_map}}\\
\subfloat[Load for the 6-zone 12-week case, all timesteps]
{\includegraphics[width=0.75\textwidth]{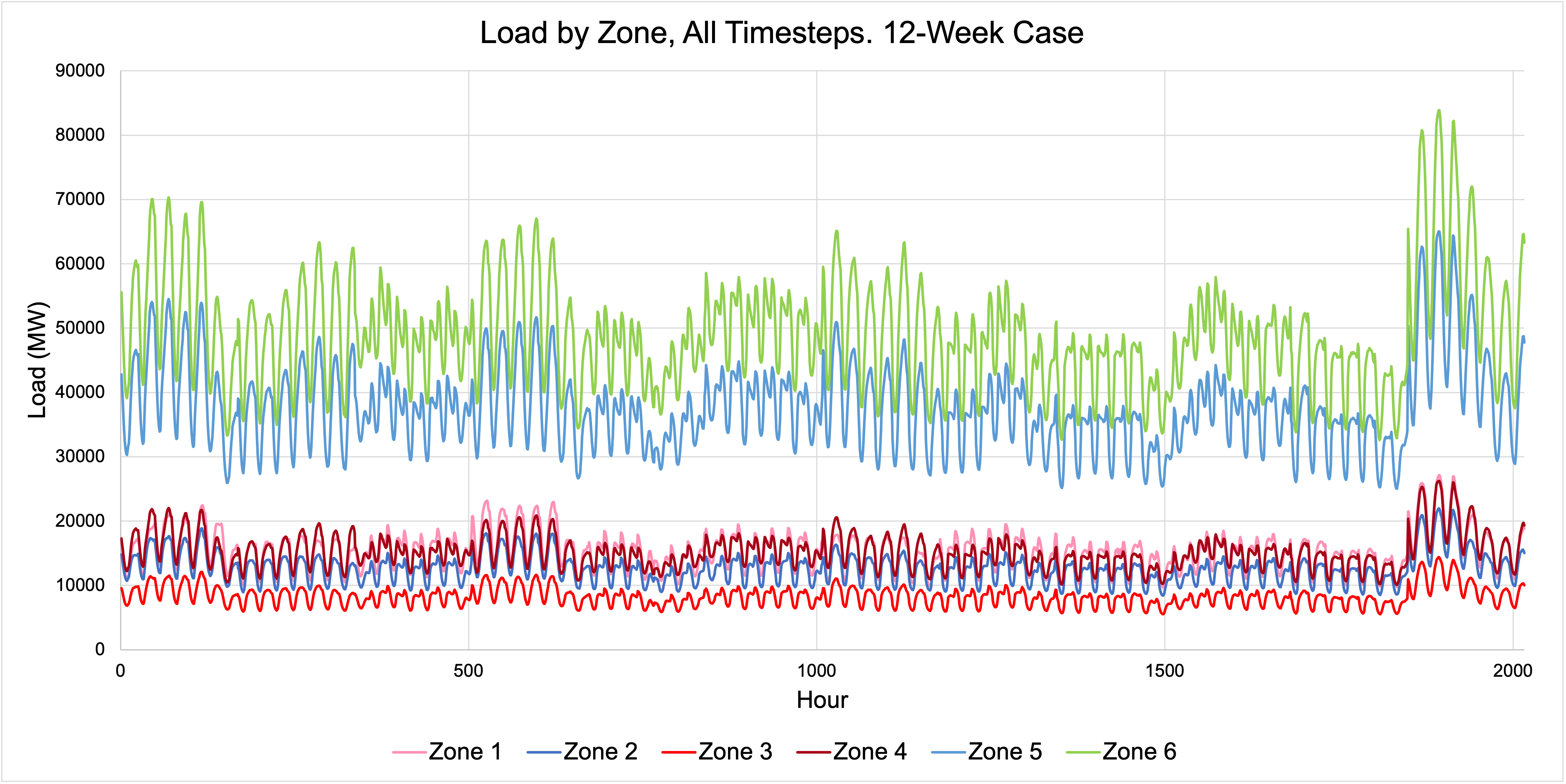}\label{load_all_timesteps}}
\end{figure}

\begin{figure}[h!]
\centering
\caption{Availability of VRE for the 6-zone 12-week case. Figure shows the first timestep of the simulation. Zones correspond to the map in figure \ref{6z_map}.}
\subfloat[Availability of solar for existing photovoltaic capacity, by zone (see \ref{6z_map}). First timestep, 6-zone 12-week case.]
{\includegraphics[width=0.4\textwidth]{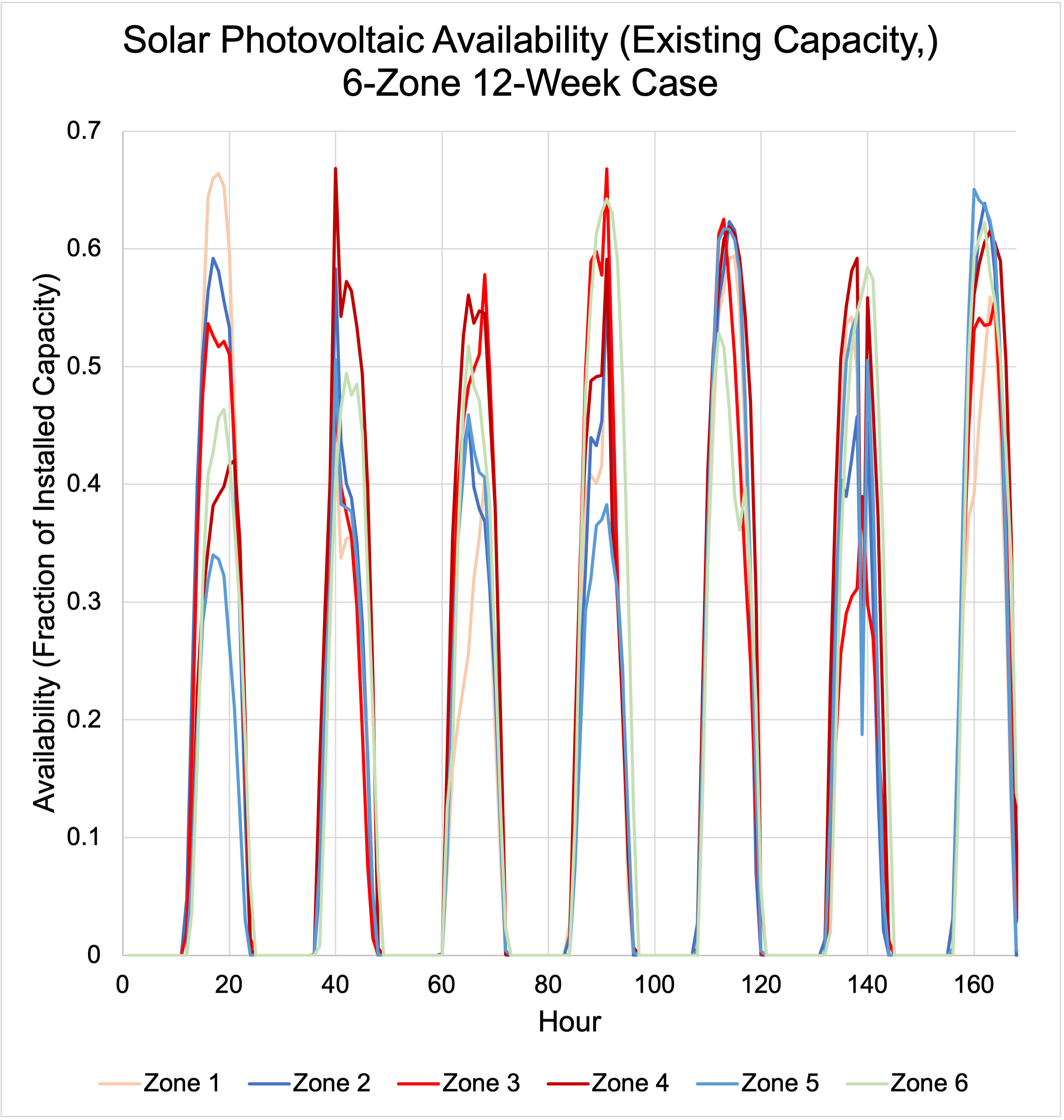}\label{solar_avail}}
\qquad
\subfloat[Availability of wind for existing wind capacity, by zone (see \ref{6z_map}). First timestep, 6-zone 12-week case.]{\includegraphics[width=0.4\textwidth]{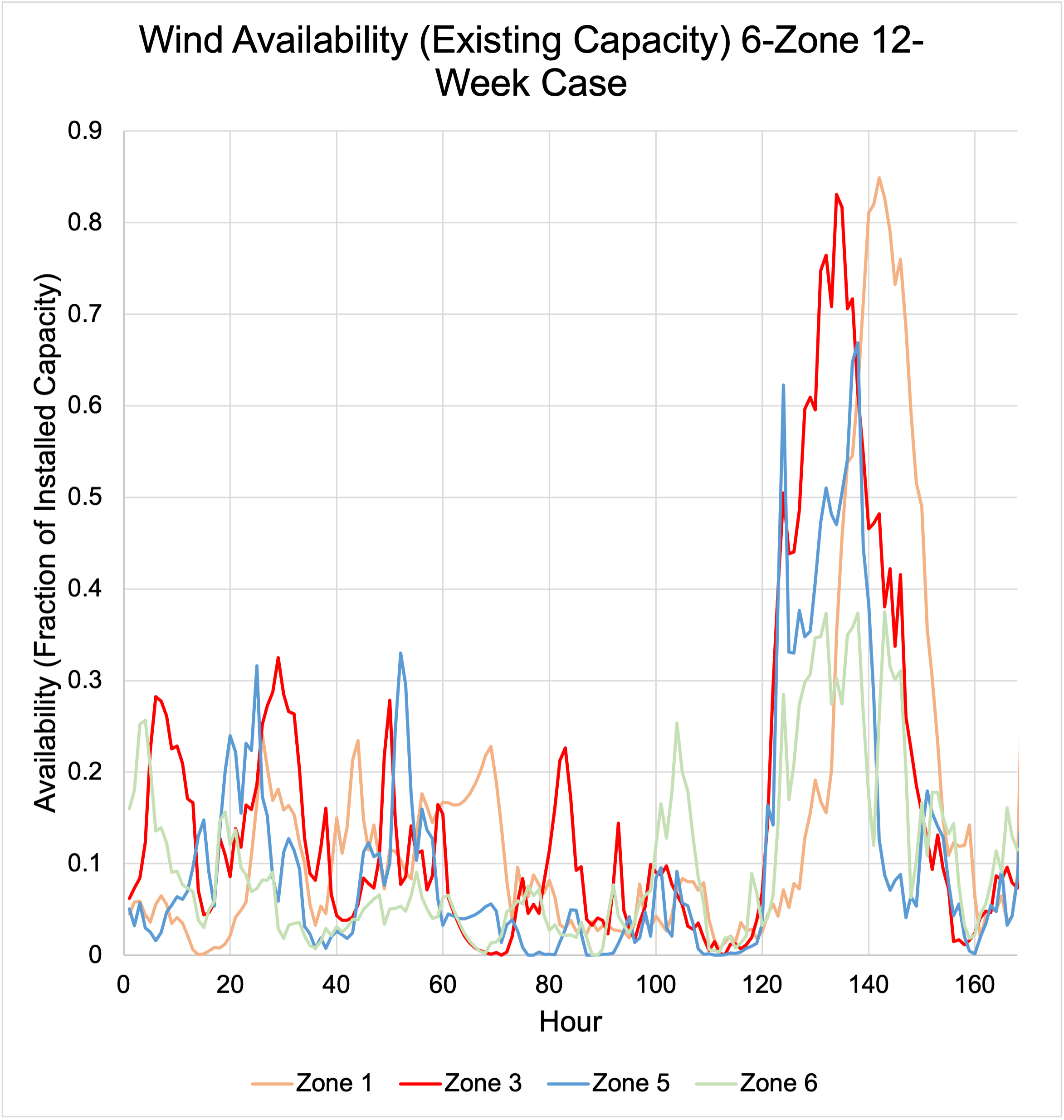}\label{wind_avail}}
\end{figure}
\cleardoublepage
\section{Results} \label{appendix_tables}
We have included the LP relaxation gap in table~\ref{lp_gap_table}. As noted on page \pageref{MSE_explanation} of the main text, trends in objective do not capture the the full implications of abstractions on models whose primary purpose is in output of variables' values and not simple system cost. As such, LP relaxation gap only tells a partial story of the impact of the transformation from LP to MILP formulations. 

\cite{donohoo2014design} notes that in linearizing investment decisions in transmission, costs for expansion can only be considered on a per-MW basis. Models will always select the from the pool of the cheapest class of transmission line first, potentially inducing a bias towards smaller lines~\citep{donohoo2014design}. This may prevent full capturing of economies of scale as an emergent property of systems. Further exploration of the specific per-generator impact of linearization on investment decisions is warranted but remains outside the scope of this manuscript.

\begin{table}[h!]
 \renewcommand{\arraystretch}{1.45}
    \centering
    \begin{tabular}{C{0.04\textwidth}|C{0.15\textwidth}|C{0.15\textwidth}|C{0.15\textwidth}}
        \toprule
         $\vert Z \vert$ & \multicolumn{3}{c}{LP Relaxation Gap (\%)} \\
         \cmidrule(lr){1-1} \cmidrule(lr){2-4}
         & REF & RPS & CO2 \\
         \cmidrule(lr){1-1} \cmidrule(lr){2-2}\cmidrule(lr){3-3}\cmidrule(lr){4-4}
         2 & 6.81e-2 & 7.81e-2 & 7.73e-2 \\
         6 & 8.81e-2 & 7.97e-2 & 6.65e-2 \\
         12 & 9.61e-2 & 0.118 & 0.106 \\
         19 & 9.34e-2 & 0.114 & 0.115 \\
         \bottomrule
    \end{tabular}
    \caption{LP relaxation gap (\%), formulated as gap in objective value between the LP and MILP case for each of the trials run. 100 $\cdot$ (UB\textsubscript{MILP} - LB\textsubscript{LP})/LB\textsubscript{LP}.}
    \label{lp_gap_table}
\end{table}
\newpage
\begin{sidewaystable}[h]
\centering
\small
\caption{Runtime for Benders and monolithic models (100s,) followed by ratio between monolithic and decomposition solution approaches. $\infty$ denotes a case with an intractable monolithic model. Cases that are intractable due to memory are noted with the superscript $M$. Cases that are intractable due to insufficient time are noted with a superscript $T$. Cases where the model outperformed its analogous model formulation are shown for the runtime rows, cases where the decomposed model outperformed monolithic are bolded for the ratio cases. Results shown for the reference case. Additional cases are included in the appendix (Table \ref{runtime_table_rps}) and in the main text (Table \ref{runtime_table_co2}, section \ref{results}.) \label{runtime_table_ref}}
 \renewcommand{\arraystretch}{1.45}
 \begin{tabular}{c|r|r|P{0.035\textheight}P{0.035\textheight}P{0.035\textheight}P{0.035\textheight}P{0.035\textheight}P{0.035\textheight}|P{0.035\textheight}P{0.035\textheight}P{0.035\textheight}P{0.035\textheight}P{0.035\textheight}P{0.035\textheight}}
\toprule
\multicolumn{1}{c}{} & \multicolumn{1}{c}{$\vert Z \vert$} & \multicolumn{1}{c}{$\vert G \vert$} & \multicolumn{6}{c}{LP} & \multicolumn{6}{c}{MILP} \\
\cmidrule(lr){1-3} \cmidrule(lr){4-9} \cmidrule(lr){10-15}
\multicolumn{3}{r|}{Weeks $\rightarrow$} & 2 & 12 & 22 & 32 & 42 & 52 & 2 & 12 & 22 & 32 & 42 & 52 \\
\midrule
 & 2 & 62 & 0.9 & 1.0 & \textbf{1.2} & \textbf{1.3} & \textbf{1.5} & \textbf{1.8} & 0.9 & \textbf{1.0} & \textbf{1.2} & \textbf{1.3} & \textbf{1.5} & \textbf{1.6} \\
 & 6 & 175 & 3.4 & 6.9 & 10.7 & \textbf{15.7} & \textbf{20.7} & \textbf{22.9} & 4.1 & \textbf{8.0} & \textbf{11.8} & \textbf{17.2} & \textbf{22.4} & \textbf{22.3} \\
&12 & 285 & 31.3 & 68.2 & 105.4 & 134.9 & 175.9 & 192.0 & 34.6 & \textbf{89.0} & \textbf{115.9} & \textbf{158.0} & \textbf{210.0} & \textbf{196.4} \\
\parbox[t]{2mm}{\multirow{-4}{*}{\rotatebox[origin=c]{90}{Benders (100s)}}}$\; \;$&19 & 437 & 231.0 & 408.9 & 522.0 & 748.4 & 772.1 & 900.7 & 194.1 & \textbf{439.9} & \textbf{866.1} & \textbf{1107.5} & \textbf{1336.8} & \textbf{1235.9} \\
\midrule
& 2 & 62 & \textbf{0.4} & \textbf{0.9} & 1.7 & 3.0 & 3.7 & 5.8 & \textbf{0.6} & 5.6 & 6.5 & 13.7 & 28.4 & 41.1 \\
& 6 & 175 & \textbf{0.6} & \textbf{4.4} & \textbf{9.5} & 18.4 & 26.8 & 43.1 & \textbf{1.3} & 52.1 & 165.2 & 188.1 & 437.0 & 836.7 \\
& 12 & 285 & \textbf{1.0} & \textbf{9.9} & \textbf{29.0} & 47.7 & 75.3 & 102.8 & \textbf{3.9} & 154.4 & 950.0 & 1591.5 & 1888.6 & $\infty^T$ \\
\parbox[t]{2mm}{\multirow{-4}{*}{\rotatebox[origin=c]{90}{Mono. (100s)}}}$\; \;$&19 & 437 & \textbf{1.5} & \textbf{20.4} & \textbf{54.8} & \textbf{104.7} & \textbf{160.1} & \textbf{222.2} & \textbf{11.2} & 660.0 & $\infty^T$ & $\infty^T$ & $\infty^T$ & $\infty^T$ \\
\midrule
& 2 & 62 & 0.4 & 0.9 & \textbf{1.4} & \textbf{2.3} & \textbf{2.5} & \textbf{3.2} & 0.6 & \textbf{5.3} & \textbf{5.6} & \textbf{10.3} & \textbf{18.7} & \textbf{25.5} \\
& 6 & 175 & 0.2 & 0.6 & 0.9 & \textbf{1.2} & \textbf{1.3} & \textbf{1.9} & 0.3 & \textbf{6.6} & \textbf{14.0} & \textbf{10.9} & \textbf{19.4} & \textbf{37.4} \\
&12 & 285 & $<$0.1 & 0.1 & 0.3 & 0.3 & 0.4 & 0.6 & 0.1 & \textbf{1.7} & \textbf{8.2} & \textbf{10.1} & \textbf{9.0} & $\bm{\infty^T}$ \\
\parbox[t]{2mm}{\multirow{-4}{*}{\rotatebox[origin=c]{90}{Ratio}}}$\; \;$&19 & 437 & $<$0.1 & 0.1 & 0.1 & 0.1 & 0.2 & 0.3 & 0.1 & \textbf{1.5} & $\bm{\infty^T}$ & $\bm{\infty^T}$ & $\bm{\infty^T}$ & $\bm{\infty^T}$ \\
\bottomrule
\end{tabular}
\end{sidewaystable}
\newpage
\begin{sidewaystable}[h!]
\centering
\small
\caption{Runtime for Benders and monolithic models (100s,) followed by ratio between monolithic and decomposition solution approaches. $\infty$ denotes a case with an intractable monolithic model. Cases that are intractable due to memory are noted with the superscript $M$. Cases that are intractable due to insufficient time are noted with a superscript $T$. Cases where the model outperformed its analogous model formulation are shown for the runtime rows, cases where the decomposed model outperformed monolithic are bolded for the ratio cases. Results shown for the RPS-constrained case. Additional cases are included in the appendix (Table \ref{runtime_table_ref}) and in the main text (Table \ref{runtime_table_co2}, section \ref{results}.) \label{runtime_table_rps}}
 \renewcommand{\arraystretch}{1.45}
 \begin{tabular}{c|r|r|P{0.035\textheight}P{0.035\textheight}P{0.035\textheight}P{0.035\textheight}P{0.035\textheight}P{0.035\textheight}|P{0.035\textheight}P{0.035\textheight}P{0.035\textheight}P{0.035\textheight}P{0.035\textheight}P{0.035\textheight}}
\toprule
\multicolumn{1}{c}{} & \multicolumn{1}{c}{$\vert Z \vert$} & \multicolumn{1}{c}{$\vert G \vert$} & \multicolumn{6}{c}{LP} & \multicolumn{6}{c}{MILP} \\
\cmidrule(lr){1-3} \cmidrule(lr){4-9} \cmidrule(lr){10-15}
\multicolumn{3}{r|}{Weeks $\rightarrow$} & 2 & 12 & 22 & 32 & 42 & 52 & 2 & 12 & 22 & 32 & 42 & 52 \\
\midrule
& 2 & 62 & 0.9 & 1.1 & \textbf{1.4} & \textbf{1.6} & \textbf{1.9} & \textbf{2.3} & 1.0 & \textbf{1.3} & \textbf{1.5} & \textbf{1.9} & \textbf{2.2} & \textbf{2.2} \\
&6 & 175 & 3.7 & 8.3 & 13.8 & \textbf{16.2} & \textbf{20.0} & \textbf{24.3} & 3.6 & \textbf{10.1} & \textbf{16.3} & \textbf{22.2} & \textbf{28.6} & \textbf{38.3} \\
&12 & 285 & 33.1 & 61.2 & 1.2e2 & 1.3e2 & 1.6e2 & 1.8e2 & \textbf{36.2} & \textbf{82.9} & \textbf{1.4e2} & \textbf{1.8e2} & \textbf{2.2e2} & \textbf{2.2e2} \\
\parbox[t]{2mm}{\multirow{-4}{*}{\rotatebox[origin=c]{90}{Benders (100s)}}}$\; \;$&19 & 437 & 1.6e2 & 3.8e2 & 5.1e2 & 6.3e2 & 8.7e2 & 9.4e2 & 1.7e2 & \textbf{4.5e2} & \textbf{6.0e2} & \textbf{1.2e3} & \textbf{1.3e3} & \textbf{1.2e3} \\
\midrule
& 2 & 62 & \textbf{0.4} & \textbf{1.0} & 1.7 & 2.8 & 3.8 & 5.1 & \textbf{0.6} & 6.9 & 10.2 & 20.6 & 44.4 & 2.1e2  \\
&6 & 175 & \textbf{0.6} & \textbf{3.7} & \textbf{9.3} & 17.4 & 25.2 & 32.1 & \textbf{2.2} & 28.1 & $\infty^T$ & 3.1e2 & 4.4e2 & 8.7e2 \\
\cellcolor{white}&12 & 285 & \textbf{1.0} & \textbf{9.7} & \textbf{24.8} & \textbf{47.5} & \textbf{76.3} & \textbf{1.1e2} & 37.9 & $\infty^T$ & $\infty^T$ & $\infty^T$ & $\infty^T$ & $\infty^T$ \\
\parbox[t]{2mm}{\multirow{-4}{*}{\cellcolor{white}\rotatebox[origin=c]{90}{Mono. (100s)}}}$\; \;$&19 & 437 & \textbf{1.7} & \textbf{21.8} & \textbf{55.3} & \textbf{1.0e2} & \textbf{1.6e2} & \textbf{2.3e2} & \textbf{9.8} & $\infty^T$ & $\infty^T$ & $\infty^T$ & $\infty^T$ & $\infty^T$ \\
\midrule
& 2 & 62 & 0.5 & 0.9 & \textbf{1.3} & \textbf{1.8} & \textbf{2.0} & \textbf{2.3} & 0.6 & \textbf{5.5} & \textbf{6.7} & \textbf{11.1} & \textbf{20.5} & \textbf{97.3} \\
&6 & 175 & 0.2 & 0.4 & 0.7 & \textbf{1.1} & \textbf{1.3} & \textbf{1.3} & 0.6 & \textbf{2.8} & $\bm{\infty^T}$ & \textbf{14.2} & \textbf{15.3} & \textbf{22.7} \\
&12 & 285 & $<$0.1 & 0.2 & 0.2 & 0.4 & 0.5 & 0.6 & 1.0 & $\bm{\infty^T}$ & $\bm{\infty^T}$ & $\bm{\infty^T}$ & $\bm{\infty^T}$ & $\bm{\infty^T}$ \\
\parbox[t]{2mm}{\multirow{-4}{*}{\cellcolor{white}\rotatebox[origin=c]{90}{Ratio}}}$\; \;$&19 & 437 & 0.1 & 0.1 & 0.1 & 0.2 & 0.2 & 0.2 & $<$0.1 & $\bm{\infty^T}$ & $\bm{\infty^T}$ & $\bm{\infty^T}$ & $\bm{\infty^T}$ & $\bm{\infty^T}$ \\
\bottomrule
\end{tabular}
\end{sidewaystable}
\cleardoublepage
\newpage
\noindent \textbf{Author Biographies.} A biologist by training, A. Jacobson was introduced to J. Jenkins and energy systems work after a happenstance meeting at a local dog park. This paper is part of the dissertation work of A. Jacobson and reflects efforts to improve the practical performance of decomposition methods for large-scale macro-energy systems planning problems with detailed operational decisions and time coupling constraints. The work builds on and improves the performance of the nested decomposition introduced by N. Sepulveda in his dissertation (Sepulveda 2020). A. Jacobson implemented the decomposition method and model code and conducted numerical experiments. Q. Xu, J. Jenkins and N. Sepulveda contributed to the decomposition method. A. Jacobson, F. Pecci, Q. Xu, and J. Jenkins analyzed results and wrote the manuscript.\\ \\
\noindent \textbf{Anna Jacobson} is a PhD candidate at Princeton University studying mathematical modeling and energy policy. After earning her B.S. in Computer Science and Biology at Tufts University, she transitioned to the energy sector due to an interest in climate-relevant models' applications.\\
\textbf{Filippo Pecci} is an Associate Research Scholar at the Andlinger Center for Energy and the Environment at Princeton University. He is a member of Princeton ZERO Lab, working on optimization-based macro-energy systems models. Before joining ZERO Lab, Filippo earned a PhD from Imperial College London, where he was also a postdoctoral research associate. \\
\textbf{Nestor Sepulveda} is a management consultant working in strategy, technology development, investing and analytics for the energy transition. Nestor earned a PhD from the Massachusetts Institute of Technology developing methods combining operations research and analytics to guide the energy transition and cleantech development. He has authored peer-reviewed papers in Nature Energy, Joule, and Applied Energy among others and has had work featured in the New York Times, the Economist, Wall Street Journal, and Bloomberg, and more. He received an M.S. working on energy policy and economics and an M.S. in Nuclear Science and Engineering, both from MIT. Nestor served for more than 10 years as a Naval Officer with the Chilean Navy.\\
\textbf{Qingyu Xu} is a research scientist at Energy Internet Research Institute at Tsinghua University, working on long-term power system planning and China's electricity market design. He received his B.S. at Sun Yat-Sen University in 2013. He received his M.S. (2015) and PhD (2020) at Johns Hopkins University, and was a Postdoc with Princeton University's ZERO lab from 2020-2022.\\
\textbf{Jesse Jenkins} is an assistant professor at Princeton University in the Department of Mechanical and Aerospace Engineering and the Andlinger Center for Energy and the Environment. He is a macro-scale energy systems engineer with a focus on the rapidly evolving electricity sector and leads the Princeton ZERO Lab, which focuses on improving and applying optimization-based energy systems models to evaluate low-carbon energy technologies and generate insights to guide policy and planning decisions. Jesse earned a PhD and SM from the Massachusetts Institute of Technology and was previously a postdoctoral environmental fellow at Harvard University.

\newpage

\bibliographystyle{informs2014} 
\bibliography{bibliography} 

\begin{thebibliography}{38}
\providecommand{\natexlab}[1]{#1}
\providecommand{\url}[1]{\texttt{#1}}
\providecommand{\urlprefix}{URL }

\bibitem[{An(2020)}]{an2020}
An K (2020) Battery electric bus infrastructure planning under demand
  uncertainty. \emph{Transportation Research Part C: Emerging Technologies}
  111:572--587, ISSN 0968090X,
  \urlprefix\url{http://dx.doi.org/10.1016/j.trc.2020.01.009}.

\bibitem[{Bezanson et~al.(2017)Bezanson, Edelman, Karpinski, \protect\BIBand{}
  Shah}]{julia_2017}
Bezanson J, Edelman A, Karpinski S, Shah VB (2017) Julia: {A} fresh approach to
  numerical computing. \emph{SIAM Review} 59(1):65--98, ISSN 00361445,
  \urlprefix\url{http://dx.doi.org/10.1137/141000671}, arXiv: 1411.1607.

\bibitem[{Bistline et~al.(2022)Bistline, Abhyankar, Blanford, Clarke, Fakhry,
  McJeon, Reilly, Roney, Wilson, Yuan et~al.}]{bistline2022actions}
Bistline J, Abhyankar N, Blanford G, Clarke L, Fakhry R, McJeon H, Reilly J,
  Roney C, Wilson T, Yuan M, et~al. (2022) Actions for reducing us emissions at
  least 50\% by 2030. \emph{Science} 376(6596):922--924.

\bibitem[{Brown et~al.(2018)Brown, Schlachtberger, Kies, Schramm,
  \protect\BIBand{} Greiner}]{brown2018}
Brown T, Schlachtberger D, Kies A, Schramm S, Greiner M (2018) Synergies of
  sector coupling and transmission reinforcement in a cost-optimised, highly
  renewable {European} energy system. \emph{Energy} 160:720--739, ISSN
  0360-5442, \urlprefix\url{http://dx.doi.org/10.1016/j.energy.2018.06.222}.

\bibitem[{Cho et~al.(2022)Cho, Li, \protect\BIBand{} Grossmann}]{cho2022}
Cho S, Li C, Grossmann IE (2022) Recent advances and challenges in optimization
  models for expansion planning of power systems and reliability optimization.
  \emph{Computers \& Chemical Engineering} 165:107924, ISSN 00981354,
  \urlprefix\url{http://dx.doi.org/10.1016/j.compchemeng.2022.107924}.

\bibitem[{Donohoo-Vallett(2014)}]{donohoo2014design}
Donohoo-Vallett PE (2014) \emph{Design of wide-area electric transmission
  networks under uncertainty: Methods for dimensionality reduction}. Ph.D.
  thesis, Massachusetts Institute of Technology.

\bibitem[{Dunning et~al.(2017)Dunning, Huchette, \protect\BIBand{}
  Lubin}]{jump_2017}
Dunning I, Huchette J, Lubin M (2017) {JuMP}: {A} modeling language for
  mathematical optimization. \emph{SIAM Review} 59(2):295--320, ISSN 00361445,
  \urlprefix\url{http://dx.doi.org/10.1137/15M1020575}, arXiv: 1508.01982.

\bibitem[{Frew \protect\BIBand{} Jacobson(2016)}]{frew2016temporal}
Frew BA, Jacobson MZ (2016) Temporal and spatial tradeoffs in power system
  modeling with assumptions about storage: An application of the power model.
  \emph{Energy} 117:198--213.

\bibitem[{Frysztacki et~al.(2022)Frysztacki, Recht, \protect\BIBand{}
  Brown}]{frysztacki2022}
Frysztacki MM, Recht G, Brown T (2022) A comparison of clustering methods for
  the spatial reduction of renewable electricity optimisation models of
  {Europe}. \emph{Energy Informatics} 5(1):4, ISSN 2520-8942,
  \urlprefix\url{http://dx.doi.org/10.1186/s42162-022-00187-7}.

\bibitem[{Han et~al.(2021)Han, Xiong, Zeng, Feng, Zhao, \protect\BIBand{}
  Yan}]{han2021comprehensive}
Han P, Xiong X, Zeng D, Feng Y, Zhao L, Yan Z (2021) Comprehensive evaluation
  of various linear programming solvers for sced in power systems. \emph{2021
  IEEE 5th Conference on Energy Internet and Energy System Integration (EI2)},
  1058--1062 (IEEE).

\bibitem[{He et~al.(2021)He, Mallapragada, Bose, Heuberger-Austin,
  \protect\BIBand{} Gençer}]{he2021}
He G, Mallapragada DS, Bose A, Heuberger-Austin CF, Gençer E (2021) Sector
  coupling \textit{via} hydrogen to lower the cost of energy system
  decarbonization. \emph{Energy \& Environmental Science} 14(9):4635--4646,
  ISSN 1754-5692, 1754-5706,
  \urlprefix\url{http://dx.doi.org/10.1039/D1EE00627D}.

\bibitem[{Helist{\"o} et~al.(2021)Helist{\"o}, Kiviluoma, Morales-Espa{\~n}a,
  \protect\BIBand{} O’Dwyer}]{helisto2021impact}
Helist{\"o} N, Kiviluoma J, Morales-Espa{\~n}a G, O’Dwyer C (2021) Impact of
  operational details and temporal representations on investment planning in
  energy systems dominated by wind and solar. \emph{Applied Energy} 290:116712.

\bibitem[{Lara et~al.(2018)Lara, Mallapragada, Papageorgiou, Venkatesh,
  \protect\BIBand{} Grossmann}]{lara2018deterministic}
Lara CL, Mallapragada DS, Papageorgiou DJ, Venkatesh A, Grossmann IE (2018)
  Deterministic electric power infrastructure planning: Mixed-integer
  programming model and nested decomposition algorithm. \emph{European Journal
  of Operational Research} 271(3):1037--1054.

\bibitem[{Larson et~al.(2021)Larson, Greig, Jenkins, Mayfield, Pascale, Zhang,
  Drossman, Williams, Pacala, Socolow, Baik, Birdsey, Duke, Jones, Haley,
  Leslie, Paustian, \protect\BIBand{} Swang}]{nzap}
Larson E, Greig C, Jenkins J, Mayfield E, Pascale A, Zhang C, Drossman J,
  Williams R, Pacala S, Socolow R, Baik E, Birdsey R, Duke R, Jones R, Haley B,
  Leslie E, Paustian K, Swang A (2021) Net zero america: Potential pathways,
  infrastructure, and impacts, final report summary.

\bibitem[{Li et~al.(2022)Li, Conejo, Liu, Omell, Siirola, \protect\BIBand{}
  Grossmann}]{li2022mixed}
Li C, Conejo AJ, Liu P, Omell BP, Siirola JD, Grossmann IE (2022) Mixed-integer
  linear programming models and algorithms for generation and transmission
  expansion planning of power systems. \emph{European Journal of Operational
  Research} 297(3):1071--1082.

\bibitem[{Lohmann \protect\BIBand{} Rebennack(2017)}]{lohmann2017tailored}
Lohmann T, Rebennack S (2017) Tailored benders decomposition for a long-term
  power expansion model with short-term demand response. \emph{Management
  Science} 63(6):2027--2048.

\bibitem[{Mallapragada et~al.(2018)Mallapragada, Papageorgiou, Venkatesh, Lara,
  \protect\BIBand{} Grossmann}]{mallapragada2018}
Mallapragada DS, Papageorgiou DJ, Venkatesh A, Lara CL, Grossmann IE (2018)
  Impact of model resolution on scenario outcomes for electricity sector system
  expansion. \emph{Energy} 163:1231--1244, ISSN 0360-5442,
  \urlprefix\url{http://dx.doi.org/https://doi.org/10.1016/j.energy.2018.08.015}.

\bibitem[{Mallapragada et~al.(2020)Mallapragada, Sepulveda, \protect\BIBand{}
  Jenkins}]{mallapragada2020}
Mallapragada DS, Sepulveda NA, Jenkins JD (2020) Long-run system value of
  battery energy storage in future grids with increasing wind and solar
  generation. \emph{Applied Energy} 275:115390, ISSN 0306-2619,
  \urlprefix\url{http://dx.doi.org/https://doi.org/10.1016/j.apenergy.2020.115390}.

\bibitem[{Mittelmann(2023)}]{mittelmann2023}
Mittelmann H (2023) Benchmarks for optimization software.
  \urlprefix\url{http://plato.asu.edu/bench.html}, accessed on 5.1.2023.

\bibitem[{Munoz et~al.(2016)Munoz, Hobbs, \protect\BIBand{}
  Watson}]{munoz2016new}
Munoz FD, Hobbs BF, Watson JP (2016) New bounding and decomposition approaches
  for milp investment problems: Multi-area transmission and generation planning
  under policy constraints. \emph{European Journal of Operational Research}
  248(3):888--898.

\bibitem[{Naderi \protect\BIBand{} Pishvaee(2017)}]{naderi2017}
Naderi MJ, Pishvaee MS (2017) A stochastic programming approach to integrated
  water supply and wastewater collection network design problem.
  \emph{Computers \& Chemical Engineering} 104:107--127, ISSN 00981354,
  \urlprefix\url{http://dx.doi.org/10.1016/j.compchemeng.2017.04.003}.

\bibitem[{Neumann et~al.(2022)Neumann, Hagenmeyer, \protect\BIBand{}
  Brown}]{neumann2022}
Neumann F, Hagenmeyer V, Brown T (2022) Assessments of linear power flow and
  transmission loss approximations in coordinated capacity expansion problems.
  \emph{Applied Energy} 314:118859, ISSN 03062619,
  \urlprefix\url{http://dx.doi.org/10.1016/j.apenergy.2022.118859}.

\bibitem[{Palmintier \protect\BIBand{}
  Webster(2013)}]{palmintier2013heterogeneous}
Palmintier BS, Webster MD (2013) Heterogeneous unit clustering for efficient
  operational flexibility modeling. \emph{IEEE Transactions on Power Systems}
  29(3):1089--1098.

\bibitem[{Palmintier \protect\BIBand{} Webster(2015)}]{palmintier2015impact}
Palmintier BS, Webster MD (2015) Impact of operational flexibility on
  electricity generation planning with renewable and carbon targets. \emph{IEEE
  Transactions on Sustainable Energy} 7(2):672--684.

\bibitem[{Patankar \protect\BIBand{} Jenkins(2020)}]{patankar2020land}
Patankar N, Jenkins J (2020) Land use, transmission expansion, and supply-chain
  scale-up implications of alternative scenarios for 100\% carbon-free
  electricity generation in the american west. \emph{AGU Fall Meeting
  Abstracts}, volume 2020, GC080--08.

\bibitem[{Pfenninger(2017)}]{pfenninger2017dealing}
Pfenninger S (2017) Dealing with multiple decades of hourly wind and pv time
  series in energy models: A comparison of methods to reduce time resolution
  and the planning implications of inter-annual variability. \emph{Applied
  energy} 197:1--13.

\bibitem[{Poncelet et~al.(2020)Poncelet, Delarue, \protect\BIBand{}
  D’haeseleer}]{poncelet2020unit}
Poncelet K, Delarue E, D’haeseleer W (2020) Unit commitment constraints in
  long-term planning models: Relevance, pitfalls and the role of assumptions on
  flexibility. \emph{Applied Energy} 258:113843.

\bibitem[{Poncelet et~al.(2016)Poncelet, Delarue, Six, Duerinck,
  \protect\BIBand{} D’haeseleer}]{poncelet2016impact}
Poncelet K, Delarue E, Six D, Duerinck J, D’haeseleer W (2016) Impact of the
  level of temporal and operational detail in energy-system planning models.
  \emph{Applied Energy} 162:631--643.

\bibitem[{Ricks et~al.(2022)Ricks, Norbeck, \protect\BIBand{}
  Jenkins}]{ricks2022}
Ricks W, Norbeck J, Jenkins J (2022) The value of in-reservoir energy storage
  for flexible dispatch of geothermal power. \emph{Applied Energy} 313:118807,
  ISSN 0306-2619,
  \urlprefix\url{http://dx.doi.org/10.1016/j.apenergy.2022.118807}.

\bibitem[{Ricks et~al.(2023)Ricks, Xu, \protect\BIBand{} Jenkins}]{ricks2023}
Ricks W, Xu Q, Jenkins JD (2023) Minimizing emissions from grid-based hydrogen
  production in the {United} {States}. \emph{Environmental Research Letters}
  18(1):014025, ISSN 1748-9326,
  \urlprefix\url{http://dx.doi.org/10.1088/1748-9326/acacb5}.

\bibitem[{Ringkj{\o}b et~al.(2018)Ringkj{\o}b, Haugan, \protect\BIBand{}
  Solbrekke}]{ringkjob2018review}
Ringkj{\o}b HK, Haugan PM, Solbrekke IM (2018) A review of modelling tools for
  energy and electricity systems with large shares of variable renewables.
  \emph{Renewable and Sustainable Energy Reviews} 96:440--459.

\bibitem[{Schivley et~al.(2021)Schivley, Welty, \protect\BIBand{}
  Patankar}]{schivley2021powergenome}
Schivley G, Welty E, Patankar N (2021) Powergenome/powergenome: v0. 4.1.

\bibitem[{Sepulveda(2020)}]{sepulveda2020decarbonization}
Sepulveda NA (2020) \emph{Decarbonization of power systems, multi-stage
  decision-making with policy and technology uncertainty}. Ph.D. thesis,
  Massachusetts Institute of Technology.

\bibitem[{Shah \protect\BIBand{} Ierapetritou(2012)}]{shah2012}
Shah NK, Ierapetritou MG (2012) Integrated production planning and scheduling
  optimization of multisite, multiproduct process industry. \emph{Computers \&
  Chemical Engineering} 37:214--226, ISSN 00981354,
  \urlprefix\url{http://dx.doi.org/10.1016/j.compchemeng.2011.08.007}.

\bibitem[{Siala \protect\BIBand{} Mahfouz(2019)}]{siala2019impact}
Siala K, Mahfouz MY (2019) Impact of the choice of regions on energy system
  models. \emph{Energy Strategy Reviews} 25:75--85.

\bibitem[{Victoria et~al.(2022)Victoria, Zeyen, \protect\BIBand{}
  Brown}]{victoria2022}
Victoria M, Zeyen E, Brown T (2022) Speed of technological transformations
  required in {Europe} to achieve different climate goals. \emph{Joule}
  6(5):1066--1086, ISSN 2542-4351,
  \urlprefix\url{http://dx.doi.org/10.1016/j.joule.2022.04.016}.

\bibitem[{Victoria et~al.(2020)Victoria, Zhu, Brown, Andresen,
  \protect\BIBand{} Greiner}]{victoria2020}
Victoria M, Zhu K, Brown T, Andresen GB, Greiner M (2020) The role of
  photovoltaics in a sustainable {European} energy system under variable {CO2}
  emissions targets, transmission capacities, and costs assumptions.
  \emph{Progress in Photovoltaics: Research and Applications} 28(6):483--492,
  \urlprefix\url{http://dx.doi.org/10.1002/pip.3198}.

\bibitem[{Xu \protect\BIBand{} Hobbs(2019)}]{xu2019value}
Xu Q, Hobbs BF (2019) Value of model enhancements: quantifying the benefit of
  improved transmission planning models. \emph{IET Generation, Transmission \&
  Distribution} 13(13):2836--2845.

\end{thebibliography}

\end{document}